\documentclass[a4paper, 10pt]{amsart}
\usepackage[utf8]{inputenc}
\usepackage{amsmath, amssymb, hyperref}
\usepackage{tikz-cd}
\hypersetup{hidelinks}

\theoremstyle{plain}
\newtheorem{theorem}{\bf Theorem}[section]

\newtheorem{lemma}[theorem]{\bf Lemma}

\theoremstyle{definition}

\newcommand{\N}{\mathbb N}
\newcommand{\Z}{\mathbb Z}
\newcommand{\R}{\mathbb R}

\newcommand{\F}{\mathbb F}
\newcommand{\C}{\mathbb C}

 \DeclareMathOperator{\ord}{ord}
\DeclareMathOperator{\lcm}{lcm} 
 
 \DeclareMathOperator{\supp}{supp}
 
\DeclareMathOperator{\Hom}{Hom}
 
\DeclareMathOperator{\im}{im}

\numberwithin{equation}{section}

\date{}
\makeatletter
\@ifundefined{subjclassname@2020}{%
 \@namedef{subjclassname@2020}{\textup{2020} Mathematics Subject Classification}}{}
\makeatother

\begin{document}

\title[Extremal minimal zero-sum sequences]{On the order of elements \\ in extremal minimal zero-sum sequences}

\author{Alfred Geroldinger \and Yongke Qu \and Guoqing Wang\textsuperscript{*}}
\thanks{*Corresponding author, was supported by the National Natural Science Foundation
of China, grant no. 12371335}

\address{Department of Mathematics and Scientific Computing\\ University of Graz, NAWI Graz\\ Heinrichstra{\ss}e 36\\ 8010 Graz, Austria}
\email{alfred.geroldinger@uni-graz.at}
\urladdr{https://geroldinger.github.io/}

\address{Department of Mathematics \\ Luoyang Normal University, Luoyang 471934, P.R. China}
\email{yongke1239@163.com}

\address{School of Mathematical Sciences \\ Tiangong University \\ Tianjin 300387, P.R. China}
\email{gqwang1979@aliyun.com}

\renewcommand{\keywordsname}{Keywords}


\begin{abstract}
Let $G = C_{n_1} \oplus \ldots \oplus C_{n_r}$ with $1 < n_1 \mid \ldots \mid n_r$ be a finite abelian group. The order of elements in minimal zero-sum sequences over $G$ of maximal length is a widely studied topic in zero-sum theory. A key focus is placed on the questions of whether there is always an element of order $n_r = \exp (G)$ and whether the orders of all elements are multiples of $n_1$. We provide further classes of groups where both questions have an affirmative answer but we also provide-for the first time-classes of groups where both questions have a negative answer. This disproves long-standing conjectures.
\end{abstract}

\maketitle

\section{Introduction}

Let $G = C_{n_1} \oplus \ldots \oplus C_{n_r}$ with $1 < n_1 \mid \ldots \mid n_r$ be a finite abelian group. The (large) Davenport constant $\mathsf D (G)$ is the maximal length of a minimal zero-sum sequence over $G$. A straightforward example shows that
$$
\mathsf D^* (G) :=  1 + \sum_{i=1}^r (n_i-1) \le \mathsf D (G) \,.
$$
Equality holds for $p$-groups, for groups of rank at most two, and others (e.g., \cite{Bh-SP07a}, or \cite{Gi18a} for an asymptotic result). There are also various series of groups with strict inequality (see \cite{Ge-Gr-Zh26a} for an overview and \cite{Li20a, Wa28a} for recent progress). It is still open whether equality holds for rank three groups.

Alongside the problem of determining the exact value of $\mathsf D(G)$, the corresponding inverse problem—which seeks the structure of minimal zero-sum sequences of length $\mathsf D(G)$—has also found wide attention. It is folklore that every minimal zero-sum sequence of maximal length over a cyclic group $G$ consists of one generating element, which is repeated $|G| = \exp (G)$ times. If $G$ has rank two, then the structure of minimal zero-sum sequences of length $\mathsf D (G)$ is fully characterized (\cite[Theorem 4.5.2]{Ge-Gr-Zh26a}). In particular, every such sequence has one element of order $\exp (G)$. Apart from some small groups and from groups of the form $G=C_2 \oplus C_2 \oplus C_{2n}$, with $n$ odd (see \cite{Sc11b}), a full characterization of all minimal zero-sum sequences of length $\mathsf D (G)$ is unknown and seems to be out of reach with present methods.

Therefore, structural issues have so far been approached mainly from the following two angles (also depending on the groups under consideration).
\begin{itemize}
\item[(a)] What is the order of elements in   minimal zero-sum sequences of length $\mathsf D (G)$?  Is there always an element with high order, say with order $\exp (G)$? Is the order of each element a multiple of $n_1$? Apart from being of interest in its own right, this question plays a crucial role in the study of cross numbers of a given group (see \cite{Ba-Sc25a, Zh28a} for  recent contributions).

\item[(b)] What is the multiplicity of elements in minimal zero-sum sequences of length $\mathsf D (G)$? Is there always an element with high multiplicity?
    \end{itemize}
It is known since more than 20 years that there are $p$-groups of high rank having squarefree minimal zero-sum sequences  of length $\mathsf D (G)$. Substantial progress in this direction was achieved by Ordaz, Philipp, Santos, and Schmid \cite{O-P-S-S11}.

We focus on Problem (a). The first question is  whether or not every minimal zero-sum sequence of maximal length resp.  whether every zero-sum free sequence of maximal length has at least one element of order $\exp (G)$. This was either formulated as an open problem resp. as a conjecture  in various sources over the last three decades (e.g., see \cite[Conjecture~6.1]{Ga-Ge99} for an early source, \cite[Problem 3]{Gi10a}, and  \cite[Appendix B, Problem 4]{Ge-Gr-Zh26a}).  Beyond the above mentioned groups (groups of rank two and of the form $G=C_2^2 \oplus C_{2n}$), the question has an affirmative answer for various series of groups, including groups of the form $C_2^4 \oplus C_{2n}$ with $n > 70$ odd (\cite[Theorem 5.9]{Sa-Ch14a})  and all $p$-groups (\cite[Proposition 3.2.14]{Ge-Gr-Zh26a}). Girard could even show that every zero-sum free sequence of maximal length  over a finite abelian $p$-group contains at least $(\exp(G)-1)$ terms of maximal order~\cite[Corollary 2.3]{Gi10a}.
Moreover, it was proved that every element of a minimal zero-sum sequence of maximal length over a $p$-group of the form $G = C_{p^k}^r$ has order $\exp (G) = p^k$ (see \cite[Theorem 5.3]{Gi10a} and
\cite[Corollary 3.6]{Ge-Ya27a}).

Schmid observed that there are minimal zero-sum sequences of length $\mathsf D^*(G)$ with no term of maximal order (\cite[discussion following Corollary~4.6]{Sc11b}). However, since in general  $\mathsf D(G)>\mathsf D^*(G)$ may hold,  the essential issue is to obtain the required property of the order of elements at the actual extremal length $\mathsf D (G)$, even when the precise value of $\mathsf D (G)$ is not known.

Theorem \ref{inv:main} provides the first series of groups $G$ having minimal zero-sum sequences of maximal length, which have no elements of order $\exp (G)$.  Indeed, for every fixed finite abelian group $G_0$ of exponent $2^ap^b$, where $a\ge2$, $b\ge1$, and $p$ is an odd prime, the group $G_r := C_2^r\oplus G_0$ shows the required phenomenon for a set of integers $r\in\N_0$ of natural density one. Note that, for large $r$, the precise value of $\mathsf D(C_2^r \oplus C_{2n})$ with odd $n\ge 3$ is not known, but the difference $\mathsf D (C_2^r \oplus C_{2n}) - \mathsf D^* (C_2^r \oplus C_{2n})$ goes to infinity when $r$ goes to infinity (\cite{Ma92}).

\begin{theorem}\label{inv:main}
Let $G_r=C_2^r\oplus G_0$ be a finite abelian group, where $r \in \N_0$ and $G_0$ is a finite abelian group of exponent $n=2^ap^b$ for some $a\ge2$, $b\ge1$ and $p$  an odd prime. The set of integers $r \in \N_0 $ admitting a minimal zero-sum sequence $S$ over $G_r$ such that
$$
|S|=\mathsf D(G_r)\quad \text{and} \quad \ord(g) \ \text{divides $n/2$ or $2^a \qquad $ for every $g \in \supp ( S)$}
$$
has natural density one. More precisely, the number of integers $r\in[0,x]$
for which no such sequence exists is
$O_{G_0}((\log(x+2))^{a+1})$ as $x\to\infty$.
\end{theorem}

We continue with a positive result. The groups $G_k$ in Theorem \ref{inv:general-counterpart} are `almost' $p$-groups, and it has been shown that several zero-sum invariants of $G_k$ behave as they do for $p$-groups (e.g., \cite[Theorem 3.6]{Ge-Wa-Ya27b}). Theorem \ref{inv:general-counterpart} shows that the groups $G_k$ show the same phenomenon as $p$-groups what concerns the order of elements in extremal minimal zero-sum sequences. It can also be seen as a counterpart to Theorem \ref{inv:main} (indeed, for $G=C_2^r$ and $n=2^ap^b$ with $a\ge2$, $b\ge1$, and $p$ an odd prime, take $2$ as the prime in Theorem~\ref{inv:general-counterpart}, $m=a$, and $k=p^b$. Since $\mathsf d(C_2^r)=r$, the theorem gives at least $2^a-r$ terms of order $\exp (G_k)=n$ whenever $r \in [0, 2^a-1]$).

\begin{theorem}\label{inv:general-counterpart}
Let $G_k=G\oplus C_{kp^m}\cong C_{n_1}\oplus\ldots\oplus C_{n_r}$ with $1<n_1\mid\ldots\mid n_r$, where $G$ is a $p$-group, $k,m\in\N$, and $\mathsf D(G)\le p^m$, and let $S$ be a minimal zero-sum sequence over $G_k$ with $|S|=\mathsf D(G_k)$.
\begin{enumerate}
\item $S$ has at least $p^m-\mathsf d(G)$ elements (counted with multiplicity) of order $\exp(G_k)$.

\item The order of every element of $S$ is a multiple of  $n_1$.
\end{enumerate}
\end{theorem}

Let $G = C_{n_1} \oplus \ldots \oplus C_{n_r}$ with $1 < n_1 \mid \ldots \mid n_r$.
Also the second question in Problem (a)-whether all orders of elements in minimal zero-sum sequences are multiples of $n_1$-pops up in the literature in several places. Girard formulated it as a conjecture, and it is also addressed in \cite[Appendix B, Problem 4]{Ge-Gr-Zh26a}.
Our final theorem considers the well-studied class of groups $G_r=C_2^r\oplus C_{2n}$ for odd $n\ge3$, whence $n_1=2$.  Roughly speaking, Theorem \ref{inv:girard-result} says that in the overwhelming number of cases the order of elements are multiples of $2$, but there are exceptions. This disproves Girard's Conjecture \cite[Conjecture~3]{Gi10a}.

\newpage
\begin{theorem}\label{inv:girard-result}
Let  $G_r=C_2^r\oplus C_{2n}$ with $r\in\N$ and $n \ge 3$ odd.
\begin{enumerate}
\item The set of integers $r\in\N$ for which every minimal zero-sum
sequence $S$ over $G_r$ of length $\mathsf D(G_r)$ satisfies
$$
 2\mid\ord(g)\qquad\text{for every }g\in\supp(S)
$$
has natural density one.

\item If $n$ is a prime power, then, for every sufficiently large $r$ (depending on $n$), every element in  every minimal zero-sum sequence over $G_r$ of length $\mathsf D(G_r)$ has even order.

\item If $n=3$, then there exists a minimal zero-sum sequence over $G_r$ of length $\mathsf D(G_r)$ having an element of odd order if and only if $r=3$.
\end{enumerate}
\end{theorem}

In Section \ref{sec:prelim}, we fix notation and discuss some preparatory lemmas. Next, Theorems \ref{inv:main}, \ref{inv:general-counterpart}, and \ref{inv:girard-result}  are each proven in a separate section.

\section{Background and some Lemmas} \label{sec:prelim}
Let $\N$ denote all positive integers, use $\N_0$ for all non-negative integers, and let $\mathbb P\subseteq \N$ denote the set of all positive primes. For $p \in \mathbb P$, we set $\F_p = \Z/ p \Z$. For $a, b \in \R$, we use $[a,b] = \{x \in \Z \colon a \le x \le b\}$ to denote the discrete interval between $a$ and $b$. For $k, \ell \in \N_0$, $\binom{k}{\ell}$ is the usual binomial coefficient with $\binom{k}{\ell}=0$ if $\ell > k$. We use the $O$-notation and the $\ll$-notation. For a nonempty subset $D \subseteq \C$ and for two maps $f, g \colon D \to \C$, we write
$$
f = O (g), \ f(s) = O ( g(s) ) \ \text{for all $s \in D$}, \ f \ll g, \ \text{and} \ f(s) \ll g (s) \ \text{for all $s \in D$}
$$
if there is some $M \in \R_{>0}$ such that $|f(s)| \le M |g(s)|$ for all $s \in D$.

For $n \in \N$, we denote by $C_n$ an additive  cyclic group of order $n$. Let $G = C_{n_1} \oplus \ldots \oplus C_{n_r}$ be a finite abelian group with $1 < n_1 \mid \ldots \mid n_r$. Then $n_r = \exp (G)$ is the exponent of $G$ and $r = \mathsf r (G)$ is the rank of $G$. For subsets $A,B \subseteq G$, the sumset and difference set are
$$
A+B=\{a+b\colon a\in A,\ b\in B\},\qquad
A-B=\{a-b\colon a\in A,\ b\in B\},
$$
and we set $-A=\{-a\colon a\in A\}$.  We say that  $A$ is covered by $k$ cosets of $H$ (for some subgroup $H \le G$) if $A\subseteq\bigcup_{j=1}^k(g_j+H)$ for some $g_j\in G$.

Sequences over $G$ are written multiplicatively and considered as elements of the free abelian monoid $\mathcal F (G)$. If
$$
S=g_1\cdot \ldots \cdot  g_\ell \in \mathcal F (G) \,,
$$
then $|S|=\ell$ is the length, $\sigma (S) = g_1+ \ldots + g_\ell$ is the sum,  and $\supp (S) = \{g_1, \ldots, g_\ell \} \subseteq G$ is the support of $S$.
For $T\in\mathcal F(G)$, we write $T\mid S$ if $T$ is a subsequence of $S$. In this case, $ST^{-1}$ denotes the unique sequence $U$ such that $S=TU$.
Furthermore,  $S$ is
\begin{itemize}
\item a {\it zero-sum sequence} if $\sigma (S)=0$,

\item a {\it minimal zero-sum sequence} if $\sigma (S)=0$ and $\sum_{i \in I}g_i \ne 0$ for all $\emptyset \ne I \subsetneq [1, \ell]$,

\item {\it zero-sum free} if $\sum_{i \in I}g_i \ne 0$ for all $\emptyset \ne I \subseteq [1, \ell]$.
\end{itemize}
The set of minimal zero-sum sequences is denoted by $\mathcal A(G)$. Then
\begin{itemize}
\item $\mathsf D (G)=\max\{|S| \colon S\in\mathcal A(G)\}$ is the {\it (large) Davenport constant}, and

\item $\mathsf d (G) = \max \{ |S| \colon S \in \mathcal F (G) \ \text{is zero-sum free} \}$ is the {\it (small) Davenport constant} of $G$.
\end{itemize}
We have $\mathsf D (G) = \mathsf d (G)+1$,
\begin{equation} \label{basic-1}
\mathsf d^* (G):=\sum_{i=1}^r (n_i-1) \le \mathsf d (G) \quad \text{and} \quad \mathsf D^* (G) = 1 + \sum_{i=1}^r (n_i-1) \le \mathsf D (G) \,.
\end{equation}
We will use without further mention that both inequalities are equalities for $p$-groups and for groups of rank at most two (for proofs see \cite[Chapter 5]{Ge-HK06a}).

\smallskip
We gather some basics on group algebras. Let $F$ be a field. The group algebra $F[G]$ is an $F$-vector space with basis $\{X^g \colon g \in G\}$. Multiplication is defined by $X^g X^h = X^{g+h}$ for all $g, h \in G$, and $X^0 = 1_{F [G]}$ is the identity element. Thus, every $f \in F[G]$ can be written uniquely in the form $f = \sum_{g \in G} c_gX^g$, where $c_g \in F$ for all $g \in G$, and we set
$$
[X^g]f = c_g \quad \text{for every $g \in G$} \,.
$$
The epimorphism
$$
\varepsilon \colon F[G] \to F, \qquad \text{defined by } \quad \varepsilon \Big( \sum_{g \in G} c_gX^g \Big) = \sum_{g \in G} c_g
$$
is called the {\it augmentation map}, and its kernel $I_G = \ker (\varepsilon )$ is the {\it augmentation ideal} of $F[G]$.
Moreover, the augmentation ideal is generated by the elements
$X^g-1$ with $g\in G$, that is,
\begin{equation}\label{eq:augmentation-ideal}
I_G = \langle X^g-1:g\in G \rangle \,.
\end{equation}
If $G$ is a $p$-group and $F$ has characteristic $p$, then, by \cite[Proposition 3.4.7]{Ge-Gr-Zh26a}, we have
\begin{equation} \label{vanishing-of-the-augmentation-ideal}
I_G^{\mathsf D (G)} = \{0\} \subseteq F[G] \,.
\end{equation}
For a sequence $S=g_1 \cdot \ldots \cdot g_\ell \in \mathcal F (G)$, we set
\begin{equation} \label{basic-notation}
\Pi_F(S) = \Pi (S) = \prod_{i=1}^\ell (1 - X^{g_i}) \in F[G] \,.
\end{equation}
If $\Pi (S)=0$, then $S$ is not zero-sum free.

If $V$ is a finite dimensional $F$-vector space, then $V^* = \Hom_F (V, F)$ denotes the {\it dual space} of linear forms.
Since we use methods from linear algebra and group algebras, we will always regard elementary $2$-groups  as vector spaces over  $\F_2$.

A polynomial
$f \in  F [X_1,\ldots,X_\ell]$ is called \emph{multilinear}
if the degree of each indeterminate  is at most $1$.
Equivalently, $f$ can be written in the form
$$
f(X_1,\ldots,X_\ell) = \sum_{I\subseteq [1,\ell]} c_I\prod_{i\in I}X_i \ \mbox{ with } \ c_I\in  F \,.
$$

For a set $\Omega$ and a subset $A\subseteq\Omega$, the
\emph{indicator function} $\mathbf1_A\colon\Omega\to F$ is defined by
$$
 \mathbf1_A(x)=
 \begin{cases}
  1_F & \text{if }x\in A,\\
  0_F & \text{if }x\in\Omega\setminus A.
 \end{cases}
$$

\begin{lemma} \label{2.1}
Let $G$ be a finite abelian group and let $S \in \mathcal A (G)$ with $|S|= \mathsf D (G)$.
\begin{enumerate}
\item If $H \subseteq G$ is a subgroup, then $\mathsf d (H)+\mathsf d (G/H)\le \mathsf d (G)$.

\item $\langle \supp (S) \rangle = G$.

\item  If $G = G_1 \oplus G_2$ and $\mathsf p_1 \colon G \to G_1$ is the projection, then $\langle \supp ( \mathsf p_1 (S) ) \rangle = G_1$.
\end{enumerate}
\end{lemma}

\begin{proof}
1. and 2. See \cite[Chapter 5]{Ge-HK06a}.

3. If $\langle \supp ( \mathsf p_1 (S) ) \rangle = G_1'$, then $S \in \mathcal F (G_1' \oplus G_2)$. Then Part 2 implies that
$G_1 \oplus G_2 = G = \langle \supp (S) \rangle  \subseteq G_1' \oplus G_2$, whence $G_1' = G_1$.
\end{proof}

\begin{lemma}\label{prelim:group-algebra}
Let  $G$ be a finite abelian $q$-group for some prime $q$.
Let $V$ be a finite-dimensional $\F_2$-vector space, let $\lambda_1,\ldots,\lambda_\ell\colon V\to \{0,1\} \subset \Z$ be maps for some $\ell \in \N$, let $h_1,\ldots,h_\ell\in G$, and let $F_G \colon V \to G$ be defined by
$F_G(x)=\sum_{i=1}^\ell \lambda_i(x) h_i$ for $x\in V$.
Then there is a multilinear polynomial $f\in\F_q[X_1,\ldots,X_\ell]$ of total degree
at most $\mathsf d(G)$ such that for every $x\in V$ we have
$$
 f(\lambda_1(x)+q\Z,\ldots,\lambda_\ell(x)+q\Z)=\mathbf1_{F_G^{-1} (0)}(x) \,, \quad \text{where} \quad \mathbf1_{F_G^{-1}(0)}\colon V\to\F_q \ \text{is the indicator function.}
$$
\end{lemma}

\begin{proof}
Let $x\in V$. Then
$$
 X^{F_G(x)}=\prod_{i=1}^\ell \bigl( 1+(X^{h_i}-1) \big(\lambda_i(x) +q\Z\big) \bigr) \in \F_q[G],
$$
because the $i$th factor is $X^{h_i}$ if $\lambda_i(x)=1$ and $1=X^0$ otherwise.
Expanding the product on the right gives
\begin{equation}\label{growth:subset-expansion}
 X^{F_G(x)}=\sum_{J\subseteq[1,\ell]}\left(\prod_{i\in J}(X^{h_i}-1)\right)\left(\prod_{i\in J} \big(\lambda_i(x)+q\Z)\right).
\end{equation}
The summand for  $J=\emptyset$ equals $1=X^0$.
Since $X^{h_i}-1\in I_G$ for each $i \in [1, \ell]$,
$\prod_{i\in J}(X^{h_i}-1)\in I_G^{|J|}$. By \eqref{vanishing-of-the-augmentation-ideal},  $I_G^{\mathsf d(G)+1}=\{0\}$, whence this  product equals zero
whenever $|J|>\mathsf d(G)$.

For each $x\in V$,  $F_G(x) \in G$,  $X^{F_G(x)} \in \F_q[G]$, whence
$$
[X^0]X^{F_G(x)}=\mathbf1_{F_G^{-1}(0)}(x)=
 \begin{cases}
  1_{\F_q} & \text{if } F_G(x)=0,\\
  0_{\F_q} & \text{if } F_G(x)\ne 0.
 \end{cases}
$$

Let us consider the coefficient of $X^0$ at both sides of \eqref{growth:subset-expansion}. On the left
we obtain $\mathbf1_{F_G^{-1}(0)}(x)$. In the sum on the right, we note that
$\prod_{i\in J}(\lambda_i(x)+q\Z) \in \F_q$, whence
$$
 \mathbf1_{F_G^{-1}(0)}(x)
 =\sum_{\substack{J\subseteq[1,\ell]\\ |J|\le\mathsf d(G)}}
 c_J\prod_{i\in J}(\lambda_i(x)+q\Z),
 \qquad \text{ where } \ c_J=[X^0]\!\left(\prod_{i\in J}(X^{h_i}-X^0)\right).
$$
Since $c_J \in \F_q$ does not depend on $x$, the polynomial
$$
 \sum_{\substack{J\subseteq[1,\ell]\\ |J|\le\mathsf d(G)}}
 c_J\prod_{i\in J}X_i \ \in\ \F_q[X_1,\ldots,X_\ell],
$$
is multilinear with total degree  at most $\mathsf d(G)$.
Substituting $X_i=\lambda_i(x)+q\Z$ for every $i$ gives
$\mathbf1_{F_G^{-1}(0)}(x)$ by the preceding identity.
\end{proof}

\begin{lemma}\label{prelim:cosets}
Let $V$ be a finite-dimensional $\F_2$-vector space and let  $U\le W\le V$ be subspaces.
Let $\ell\in\N$ and let $\lambda_1,\ldots,\lambda_\ell\colon V\to \{0,1\} \subset \Z$ be  maps.
Suppose that for each  $i\in[1,\ell]$ and each $w\in W$, either $\lambda_i(x+w)=\lambda_i(x)$ for every $x\in W$ or
$\lambda_i(x+w)=1-\lambda_i(x)$ for every $x\in W$.
If $q$ is a prime and if  a polynomial $f\in\F_q[X_1,\ldots,X_\ell]$ with total degree at most $t \in \N_0$
satisfies $f \big(\lambda_1(x)+q\Z,\ldots,\lambda_\ell(x)+q\Z \big)=\mathbf1_U(x)$ for all $x\in W$, then
$$
 2^{\dim W-\dim U}\le\sum_{j=0}^t\binom{\ell}{j}.
$$
\end{lemma}

\begin{proof}
For $x,w\in W$, we have $x\in w+U$ if and only if $x+w\in U$, whence
$$
\mathbf1_{w+U}(x) = f \big(\lambda_1(x+w)+q\Z,\ldots,\lambda_\ell(x+w)+q\Z \big) \,.
$$
For $w\in W$, we define a polynomial $f_w$ from $f$ as follows:  keep $X_i$ when translation by $w$ preserves $\lambda_i$ on $W$,
and replace $X_i$ by $1-X_i$ when it interchanges the values $0$ and $1$.
Reducing the translation identities for the maps $\lambda_i$ modulo $q$ gives the same substitutions in $\F_q$. Thus, for every $x\in W$,
$$
 f_w \big(\lambda_1(x)+q\Z,\ldots,\lambda_\ell(x)+q\Z \big)
 =f \big(\lambda_1(x+w)+q\Z,\ldots,\lambda_\ell(x+w)+q\Z \big)
 =\mathbf1_{w+U}(x).
$$
Thus the value of $f_w$ at $(\lambda_1(x)+q\Z,\ldots,\lambda_\ell(x)+q\Z)$ is
$1$ for $x\in w+U$ and $0$ for $x\in W\setminus(w+U)$.

For each coset $C=w+U\in W/U$, take its indicator
$\mathbf1_C\colon W\to\F_q$.
We claim that these $|W/U|$ functions are linearly
independent over $\F_q$.
Suppose that
$$\sum_{C\in W/U}a_C\mathbf1_C=0
$$
with $a_C\in\F_q$. For any coset $C_0\in W/U$,
choose $x\in C_0$. Since $\mathbf1_{C_0}(x)=1$ and
$\mathbf1_C(x)=0$ for all $C\ne C_0$, evaluating the
above relation at $x$ gives $a_{C_0}=0$.
As $C_0$ was arbitrary, all coefficients vanish,
establishing the claimed linear independence.

For each coset $C\in W/U$, fix a representative $w\in C$.
Expand $f_w$ into monomials and substitute
$X_i=\lambda_i(x)+q\Z$.
Since $\lambda_i(x)+q\Z$ is either $0_{\F_q}$ or
$1_{\F_q}$, every positive power of this residue class
equals the class itself.
Thus a monomial $X_1^{a_1}\cdot \ldots \cdot X_\ell^{a_\ell}$
occurring in $f_w$ evaluates to
$\prod_{i\in J}(\lambda_i(x)+q\Z)$, where
$J=\{i\in[1,\ell]:a_i>0\}$. Each $i\in J$ contributes at least $1$
to the degree, so $|J|\le a_1+\ldots+a_\ell\le t$.
For every $x\in W$, collecting the coefficients gives
$$
 \mathbf1_C(x)=\sum_{\substack{J\subseteq[1,\ell]\\ |J|\le t}}
 b_{C,J}\prod_{i\in J} \big(\lambda_i(x)+q\Z \big),
$$
where $b_{C,J}\in\F_q$ and the empty product equals $1$.
Note, that  for fixed $J$ and
$x\in W$, the product $\prod_{i\in J}(\lambda_i(x)+q\Z)$ is independent of $C$. Let $\F_q^W$ denote the vector space of all functions $W\to\F_q$
with pointwise addition and scalar multiplication. The formula for
$\mathbf1_C$ shows that every coset indicator belongs to the subspace
$$
 \operatorname{span}_{\F_q}\!
 \left\{\,x\mapsto\prod_{i\in J}(\lambda_i(x)+q\Z)
 \ \middle|\ J\subseteq[1,\ell],\ |J|\le t\,\right\}
 \subseteq\F_q^W .
$$
Observe that the dimension of this subspace over $\F_q$ is at most
$\sum_{j=0}^t\binom{\ell}{j}$, since there are $\binom{\ell}{j}$ index sets $J$ of size $j$. The coset indicators are linearly independent over $\F_q$,
so their number $|W/U|$ is at most this dimension. Since $W$ and $U$
are vector spaces over $\F_2$, the number of cosets is
$|W/U|=2^{\dim W-\dim U}$, as claimed.
\end{proof}

Next we shall describe the zero-sum subsequences of a given sequence in terms of vectors over $\F_2$.
We use this description in the proofs of Lemma~\ref{prelim:growth} and Theorem~\ref{inv:main}.

Let $G_0$ be a finite abelian group, let $\ell\in\N$ and $r\in\N_0$, and let
$S=\prod_{i=1}^\ell(g_i+a_i)\in\mathcal F(\F_2^r\oplus G_0)$, where
$g_i\in\F_2^r$ and $a_i\in G_0$ for every $i\in[1,\ell]$. Define the linear map
\begin{equation}\label{prelim:phi-definition}
 \varphi\colon\F_2^\ell\longrightarrow\F_2^r,\qquad
 \varphi(x)=\sum_{i=1}^{\ell}x_i g_i
 \quad\text{for }x=(x_1,\ldots,x_\ell)\in\F_2^\ell.
\end{equation}
Set
\begin{equation}\label{prelim:relation-space}
 V=\ker(\varphi)=\left\{x\in\F_2^\ell\colon \sum_{i=1}^\ell x_i g_i=0\right\}.
\end{equation}
Then
\begin{equation}\label{prelim:relation-dimension}
 \dim V=\ell-\dim\im(\varphi)\ge \ell-r,
\end{equation}
and $\dim V=\ell-r$ if and only if $g_1,\ldots,g_\ell$ span $\F_2^r$.

For $x\in\F_2^\ell$, consider the subsequence
$$
 S_x=\prod_{\substack{1\le i\le\ell\\ x_i=1}}(g_i+a_i) \quad \text{and observe that} \quad
$$
\begin{equation}\label{prelim:selected-binary-sum}
 \sum_{\substack{1\le i\le\ell\\ x_i=1}}g_i
 =\sum_{i=1}^{\ell}x_i g_i
 =\varphi(x).
\end{equation}
As $V=\ker (\varphi)$, this sum is zero in $\F_2^r$ if and only if $x\in V$.
For $x=(x_1,\ldots,x_\ell)\in V$, we define $\F_2$-linear maps
$\mu_1,\ldots,\mu_\ell\colon V\to\F_2$ and the map $F\colon V\to G_0$ by
\begin{equation}\label{prelim:coordinate-sum}
 \mu_i(x)=x_i\quad\text{for }i\in[1,\ell],\qquad
 F(x)=\sum_{\substack{1\le i\le\ell\\ \mu_i(x)=1}}a_i.
\end{equation}
If $x\in V$, then
the selected $g_i$ have sum zero  by \eqref{prelim:selected-binary-sum}. Since $\mu_i(x)=x_i$ for every
$i\in[1,\ell]$, \eqref{prelim:coordinate-sum} gives $F(x)$ as the
sum of the selected $a_i$, whence
$$
 \sigma(S_x)
 =\sum_{\substack{1\le i\le\ell\\ x_i=1}}g_i
  +\sum_{\substack{1\le i\le\ell\\ x_i=1}}a_i
 =0_{\F_2^r}+F(x).
$$
Thus, for every $x\in V$, $S_x$ has sum zero in $\F_2^r\oplus G_0$ if and only if
$F(x)=0$ in $G_0$.

As $x$ ranges over $\F_2^\ell$, $S_x$
runs through all subsequences of $S$.  In particular, $S_x$ is the empty subsequence $1$ of $S$ if and
only if $x=0 \in \F_2^\ell$, and $S_x=S$ if and only if
$x=\mathbf1=(1,\ldots,1) \in \F_2^\ell$. Hence
\begin{equation}\label{prelim:zero-sumfree-criterion}
 S\text{ is zero-sumfree}\quad\Longleftrightarrow\quad F^{-1}(0)=\{0\},
\end{equation}
and
\begin{equation}\label{prelim:atom-criterion}
 S\in\mathcal A(\F_2^r\oplus G_0)
 \quad\Longleftrightarrow\quad
 \mathbf1\in V\ \text{ and }\ F^{-1}(0)=\{0,\mathbf1\}.
\end{equation}

\begin{lemma} \label{prelim:growth}
Let $G_r=\F_2^r\oplus G_0$ be as in Theorem~\ref{inv:main}.
Then there exists a constant $C>0$, depending only on $G_0$, such that
$$
\exp (G_0)-1\le\mathsf d(G_0)\le\mathsf d(G_r)-r
\le C\log(r+2)
\qquad\text{for every }r\in\N_0.
$$
\end{lemma}

\begin{proof}
To recall, $G_0$ is a finite abelian group with exponent $\exp (G_0)=n=2^ap^b$ for some $a\ge2$, $b\ge1$ and $p$  an odd prime.
Thus, \eqref{basic-1} shows that $n-1 \le \mathsf d (G_0)$
and Lemma~\ref{2.1}.1 implies that  $\mathsf d(G_r)\ge r+\mathsf d(G_0)$. It remains to show the last inequality.

We set $G_0=Q\oplus P$, where $Q$ and $P$ are
the Sylow $2$- and $p$-subgroups of $G_0$, respectively, and then
$G_r=\F_2^r\oplus Q\oplus P$.
We choose a zero-sum free sequence
$S = \prod_{i=1}^\ell (g_i+q_i+p_i) \in \mathcal F (G_r)$ of length $|S|=\ell=\mathsf d(G_r)$, with $g_i \in \F_2^r$, $q_i \in Q$, and $p_i \in P$ for every $i \in [1, \ell]$.
With $a_i=q_i+p_i$, let $V$, $\mu_i$ and $F$ be as in
\eqref{prelim:relation-space} and \eqref{prelim:coordinate-sum}.
For each $i\in[1,\ell]$, we define the map $\lambda_i\colon V\to\Z$ by
$\lambda_i(x)=1$ if $x_i=1$ and $\lambda_i(x)=0$ if $x_i=0$,
where $x=(x_1,\ldots,x_\ell)\in V$. As in Lemma~\ref{prelim:group-algebra}, we define
$F_Q\colon V\to Q$ and $F_P\colon V\to P$ by
$$
 F_Q(x)=\sum_{i=1}^\ell\lambda_i(x)q_i
 \qquad\text{and}\qquad
 F_P(x)=\sum_{i=1}^\ell\lambda_i(x)p_i.
$$
For every $x\in V$, we have $F(x)=F_Q(x)+F_P(x)$ in $Q\oplus P$. Hence
\begin{equation}\label{growth:common-zero}
 F_Q^{-1}(0)\cap F_P^{-1}(0)=F^{-1}(0)=\{0\},
\end{equation}
where the last equality follows from \eqref{prelim:zero-sumfree-criterion} since $S$ is zero-sumfree.

We choose a subspace $W\le V$ that is maximal with respect to inclusion such that
\begin{equation}\label{growth:vanishing}
 F_P| W \colon W \to P  \quad \text{is the zero map} \,.
\end{equation}
Let $e_1,\ldots,e_\ell$ be the standard basis  of $\F_2^\ell$.

We claim that the sequence $T = \prod_{i=1}^\ell \big( (e_i+W) + q_i \big) \in \mathcal F \big( (\F_2^\ell/W)\oplus Q \big) $  is zero-sum free.
Assume to the contrary that $\sum_{i \in I} \big( (e_i+W) + q_i \big) = 0$ for some  $\emptyset \ne I\subseteq[1,\ell]$.
Setting $x=\sum_{i\in I}e_i$, we obtain that
$$
W+0_Q=0=\sum_{i\in I}\big((e_i+W)+q_i\big)
=(x+W)+\sum_{i\in I}q_i.
$$
Thus $x+W=W$, whence $x\in W\subseteq V$. By \eqref{growth:vanishing}, we have $F_P(x)=0_P$.
Since $x=\sum_{i\in I}e_i$, we have $x_i=1$ exactly when $i\in I$.
Thus the indices for which $\lambda_i(x)=1$ are precisely those in $I$, and the definition of $F_Q$ gives $F_Q(x)=\sum_{i\in I}q_i=0_Q$.
By \eqref{growth:common-zero}, we obtain $x=0$. This contradicts $I\ne\emptyset$.

Therefore,  it follows that

$$
\ell = |T| \le \mathsf d((\F_2^\ell/W)\oplus Q)=\mathsf d(\F_2^\ell/W)+\mathsf d(Q)=(\ell-\dim W)+\mathsf d(Q) \,,
$$
which  implies that
\begin{equation}\label{growth:dimension-W}
 \dim W\le\mathsf d(Q).
\end{equation}

We define a map
$$
f \colon V \to \F_p  \quad \text{by} \quad f(x)=\mathbf1_{F_P^{-1}(0)}(x) \in \F_p \,,
$$
and claim that
$$
\prod_{w\in W}f(x+w)=\mathbf1_W(x).
$$
We distinguish two cases. First, if $x\in W$, then $x+w\in W$ for every $w\in W$. By \eqref{growth:vanishing}, $F_P(x+w)=0_P$, so every factor $f(x+w)$ equals $1$.
Second, if $x\notin W$, then $\langle W,x\rangle=W\cup(x+W)$. This subspace strictly contains $W$ whence,  by maximality, $F_P$ is nonzero at some element of  $\langle W,x\rangle$.
By \eqref{growth:vanishing}, this element must lie in $x+W$. Thus $F_P(x+w)\ne0_P$ for some $w\in W$, and the corresponding factor $f(x+w)$ equals zero. Hence the product equals zero. This proves the claim.

We apply Lemma~\ref{prelim:group-algebra} with $q=p$, $G=P$, and $h_i=p_i$ for each $i\in[1,\ell]$, using the maps $\lambda_i\colon V\to\Z$ defined above.
With these choices, $F_G=F_P$, so the indicator in the lemma is $f=\mathbf1_{F_P^{-1}(0)}$.
Thus there is a multilinear polynomial $g\in\F_p[X_1,\ldots,X_\ell]$ of total degree at most $\mathsf d(P)$ such that
$$
 f(x) = g\bigl(\lambda_1(x)+p\Z,\ldots,\lambda_\ell(x)+p\Z\bigr) \quad \text{for every $x \in V$} \,.
$$
For $x,w\in V$, the $i$th coordinate of $x+w$ is $x_i+w_i$ in $\F_2$.
Thus translation by $w$ preserves $\lambda_i$ on $V$ if $w_i=0$ and replaces it by $1-\lambda_i$ on $V$ if $w_i=1$.
For each $w\in W$, let $g_w$ be the polynomial obtained from $g$ by keeping $X_i$ when $w_i=0$ and substituting $1-X_i$ for $X_i$ when $w_i=1$.
Reducing the preceding identities for $\lambda_i(x+w)$ modulo $p$ shows that
$$
 f(x+w)=g_w\bigl(\lambda_1(x)+p\Z,\ldots,\lambda_\ell(x)+p\Z\bigr) \quad \text{for every $x \in V$} \,.
$$
Each substitution has degree one, so the total degree of $g_w$ is at most $\mathsf d(P)$.
Since $\prod_{w\in W}f(x+w)=\mathbf1_W(x)$, the polynomial $\prod_{w\in W}g_w$ represents $\mathbf1_W$ on $V$ and has total degree at most $$t=|W|\mathsf d(P).$$
The coordinate identities above verify the translation hypothesis on $V$. Applying Lemma~\ref{prelim:cosets} with $q=p$ to the indicator $\mathbf1_W(x)$ of $W$ on $V$ therefore gives the following bound on the number of cosets of $W$ in $V$.
$$
 \frac{2^{\dim V}}{|W|}=2^{\dim V-\dim W}
 \le\sum_{j=0}^t\binom{\ell}{j}
 \le\sum_{j=0}^t\ell^j
 \le\sum_{j=0}^t\binom{t}{j}\ell^j
 =(\ell+1)^t .
$$
By \eqref{growth:dimension-W} and the definition of $t$, we have
$$
 |W|=2^{\dim W}\le 2^{\mathsf d(Q)}
 \qquad\text{and}\qquad
 t=|W|\mathsf d(P)\le 2^{\mathsf d(Q)}\mathsf d(P).
$$
Since $Q$ and $P$ are the Sylow subgroups of the fixed group $G_0$, these upper bounds are independent of $r$ and of the chosen zero-sum free sequence. Taking logarithms to base $2$ in $2^{\dim V}/|W|\le(\ell+1)^t$ and using \eqref{prelim:relation-dimension}, we obtain
$$
 \begin{aligned}
 \ell-r\le\dim V
 &\le \log_2|W|+t\log_2(\ell+1)\\
 &\le \mathsf d(Q)+2^{\mathsf d(Q)}\mathsf d(P)\log_2(\ell+1)\\
 &\le \bigl(\mathsf d(Q)+2^{\mathsf d(Q)}\mathsf d(P)\bigr)\log_2(\ell+2).
 \end{aligned}
$$
Since $\mathsf d(Q)+2^{\mathsf d(Q)}\mathsf d(P)$
depends only on $G_0$, and since
$\log_2(\ell+2)=\log(\ell+2)/\log 2$, we obtain
$\ell-r\ll_{G_0}\log(\ell+2)$.
Since $\ell\ge r$, we have $\ell\to\infty$ as $r\to\infty$. Dividing the bound for $\ell-r$ by $\ell$ gives
$$
 0\le 1-\frac{r}{\ell}\ll_{G_0}\frac{\log(\ell+2)}{\ell}\to0  \qquad\text{as }r\to\infty.
$$
Thus $r/\ell\to 1$ as $r\to\infty$, and hence
$\ell\le 2r$ for all sufficiently large $r$.
Consequently,
$\log(\ell+2)=O(\log(r+2))$.
Combining this with $\ell-r\ll_{G_0}\log(\ell+2)$
and recalling that $\ell=\mathsf d(G_r)$, we obtain
$$
 \mathsf d(G_r)-r\ll_{G_0}\log(r+2).
$$
Increasing $C$ if necessary to cover the finitely many small $r$,
we obtain the asserted bound.
\end{proof}

\section{Proof of  Theorem \ref{inv:main}}\label{inv:density}

Let $G_r=\F_2^r\oplus G_0$ be as in Theorem~\ref{inv:main}, and put
\begin{equation}\label{main:k-definition}
 k_r=\mathsf d(G_r)-r \qquad\text{for every }r\in\N_0.
\end{equation}
For every $s\in\N_0$, we have $G_{s+1}/G_s\cong C_2$.
Since $\mathsf d(C_2)=1$, Lemma~\ref{2.1}.1 gives
$\mathsf d(G_{s+1})\ge\mathsf d(G_s)+1$. Thus the sequence
$(k_s)_{s\ge0}$ is nondecreasing, since
\begin{equation}\label{main:rank-monotonicity}
 k_{s+1}=\mathsf d(G_{s+1})-(s+1)
 \ge \mathsf d(G_s)-s=k_s
 \qquad\text{for every }s\in\N_0.
\end{equation}
Lemma~\ref{prelim:growth} gives
\begin{equation}\label{main:k-logarithmic}
 n-1\le k_r\ll_{G_0}\log(r+2) \qquad\text{for every }r\in\N_0.
\end{equation}

\medskip
\begin{proof}[Proof of Theorem~\ref{inv:main}]

Let $r\in\N_0$ and let $S=\prod_{i=1}^\ell(g_i+a_i) \in \mathcal F (G_r)$, with $g_i \in \F_2^r$ and $a_i \in G_0$ for every $i \in [1, \ell]$,  be a minimal
zero-sum sequence of length $\ell=\mathsf D(G_r)=r+k_r+1$.
By Lemma~\ref{2.1}.3, $g_1,\ldots,g_\ell$ span $\F_2^r$.
With $V$, $\mu_i \colon V \to \F_2$, and $F \colon V \to G_0$ as in
\eqref{prelim:relation-space} and \eqref{prelim:coordinate-sum}, we obtain by \eqref{prelim:relation-dimension} and \eqref{prelim:atom-criterion} that
\begin{equation}\label{main:zero-set}
 \dim V=\ell-r=k_r+1
 \qquad\text{and}\qquad F^{-1}(0)=\{0,\mathbf1\}.
\end{equation}
Let $G_0=Q\oplus P$, where $Q$ and $P$ are its Sylow $2$- and
$p$-subgroups, respectively, and let $a_i=q_i+p_i$ with $q_i\in Q$ and $p_i\in P$ for every $i\in[1,\ell]$.

We first bound, in terms of $k_r$ and $Q$, the number of
nonzero $q_i$ needed to give the same map $F$.
Recall that $\dim V^*=\dim V=k_r+1$. Choose a basis $z_1,\ldots,z_{k_r+1}$ of $V^*$.
For each $i\in[1,\ell]$, let $J_i\subseteq[1,k_r+1]$ be the unique set
for which $\mu_i=\sum_{j\in J_i}z_j$ in $V^*$.
That is,
$$
 J_i=\bigl\{j\in[1,k_r+1]\colon
       \text{the coefficient of }z_j\text{ in the basis expansion of }\mu_i
       \text{ equals }1\bigr\}.
$$
There are
$M=\sum_{s=1}^a\binom{k_r+1}{s}$ nonempty subsets $I \subseteq [1,k_r+1]$
with $|I|\le a$. Write $Q^M$ for the direct sum of $M$ copies of the group $Q$,
one for each such subset $I$. For each $i\in[1,\ell]$, let $c_i\in Q^M$
be the element whose coordinate in the copy indexed by $I$ is given by
\begin{equation}\label{main:coefficient-definition}
 c_{i,I}=
 \begin{cases}
 (-2)^{|I|-1}q_i & \text{if }I\subseteq J_i,\\
 0_Q & \text{otherwise}.
 \end{cases}
\end{equation}
Fix $x\in V$ and $i\in[1,\ell]$, and put
$$
 J(x)=\{j\in[1,k_r+1]\colon z_j(x)=1\},
 \qquad m=|J_i\cap J(x)|.
$$
By the definition of $J(x)$, we have $z_j(x)=1$ for
$j\in J(x)$ and $z_j(x)=0$ for $j\notin J(x)$.
Consequently, with all sums in the following identity
taken in $\F_2$, we obtain
\begin{equation}\label{main:coordinate-parity}
\begin{aligned}
\mu_i(x) &=\sum_{j\in J_i}z_j(x)
=\sum_{j\in J_i\cap J(x)}z_j(x)
  +\sum_{j\in J_i\setminus J(x)}z_j(x) \\
&=\sum_{j\in J_i\cap J(x)}1
=
\begin{cases}
1 & \text{if }m\text{ is odd},\\
0 & \text{if }m\text{ is even}.
\end{cases}
\end{aligned}
\end{equation}
Then we have the following  integer identity
$$
\begin{aligned}
 1-(-1)^m
 &=1-(1-2)^{|J_i\cap J(x)|}
   =1-\prod_{j\in J_i\cap J(x)}(1-2)\\
 &=1-\sum_{I\subseteq J_i\cap J(x)}(-2)^{|I|}
   =-\sum_{\emptyset\ne I\subseteq J_i\cap J(x)}(-2)^{|I|}\\
 &=2\sum_{\emptyset\ne I\subseteq J_i\cap J(x)}(-2)^{|I|-1}.
\end{aligned}
$$
Here $I$ records the factors from which we choose $-2$.
The term for $I=\emptyset$ is $1$ and cancels the initial $1$.
We divide this identity by $2$ in $\Z$ and then multiply by $q_i\in Q$. This gives
$$
 \frac{1-(-1)^m}{2}\,q_i
 =\sum_{\emptyset\ne I\subseteq J_i\cap J(x)}(-2)^{|I|-1}q_i
 =\sum_{\substack{\emptyset\ne I\subseteq J(x)\\ |I|\le a}}c_{i,I},
$$
noticing that terms with $|I|>a$ vanish because $2^a$ divides $(-2)^{|I|-1}$ and $2^a q_i=0_Q$.
For the remaining subsets of $J(x)$, the definition \eqref{main:coefficient-definition} of $c_{i,I}$ gives zero when $I\not\subseteq J_i$.
The integer coefficient on the left is $1$ when $m$ is odd and $0$ when $m$ is even.
By \eqref{main:coordinate-parity}, these cases occur precisely when
$\mu_i(x)=1$ and $\mu_i(x)=0$, respectively.
Thus, for every $x\in V$ and $i\in[1,\ell]$, we have
\begin{equation}\label{main:individual-Q-coordinate}
 \sum_{\substack{\emptyset\ne I\subseteq J(x)\\ |I|\le a}}c_{i,I}
 =
 \begin{cases}
 q_i & \text{if }\mu_i(x)=1,\\
 0_Q & \text{if }\mu_i(x)=0.
 \end{cases}
\end{equation}
Summing \eqref{main:individual-Q-coordinate} over $i\in[1,\ell]$ gives
the second equality below,
$$
 F_Q(x)
  =\sum_{\substack{1\le i\le\ell\\ \mu_i(x)=1}}q_i
  =\sum_{i=1}^\ell
       \sum_{\substack{\emptyset\ne I\subseteq J(x)\\ |I|\le a}}c_{i,I}
  =\sum_{\substack{\emptyset\ne I\subseteq J(x)\\ |I|\le a}}
       \left(\sum_{i=1}^{\ell}c_{i,I}\right).
$$

Consider the sequence $C=c_1\cdot\ldots\cdot c_\ell\in\mathcal F(Q^M)$ and
choose a subset $A\subseteq[1,\ell]$ (possibly empty) such  that
$T=\prod_{i\in A}c_i$ is a zero-sum subsequence of $C$ of maximal length.
Then $CT^{-1}$ is zero-sumfree.
For $|I|>a$, the coefficient $(-2)^{|I|-1}$ is divisible by $\exp(Q)=2^a$,
so $(-2)^{|I|-1}q_i=0_Q$. Both cases in \eqref{main:coefficient-definition}
would therefore give $0_Q$ for such an $I$. Thus only the nonempty subsets $I$ with $|I|\le a$ need be considered.
Since $\sigma(T)=\sum_{i\in A}c_i=0$ in $Q^M$, taking each coordinate gives
\begin{equation}\label{main:zero-sum-coordinates}
 \sum_{i\in A}c_{i,I}=0_Q
 \qquad\text{for every }I\subseteq[1,k_r+1]\text{ with }1\le |I|\le a.
\end{equation}

We obtain $S'$ from $S$ by replacing $g_i+q_i+p_i$ by $g_i+p_i$
at each index $i\in A$ and retaining the original term
at every index outside $A$. Thus
\begin{equation}\label{main:modified-sequence}
 S'=\prod_{i\in A}(g_i+p_i)
    \prod_{i\in[1,\ell]\setminus A}(g_i+q_i+p_i)
    \ \in\mathcal F(G_r).
\end{equation}
At each index $i$, the $\F_2^r$-component of the term of $S'$ is still $g_i$. Thus $S$ and
$S'$ have the same space $V$ and the same coordinate maps $\mu_i$.
Let $F \colon V \to G_0$ denote the sum map of $S$ as in \eqref{prelim:coordinate-sum}.
Applying the same definition to $S'$ gives the map
$F'\colon V\to Q\oplus P$. Write $F'_Q\colon V\to Q$ and
$F'_P\colon V\to P$ for its $Q$- and $P$-components, respectively.
Thus,  for every $x\in V$,
$$
 F'(x)=F'_Q(x)+F'_P(x)
 =\sum_{\substack{i\in[1,\ell]\setminus A\\ \mu_i(x)=1}}q_i
  \mathbin{+}\sum_{\substack{1\le i\le\ell\\ \mu_i(x)=1}}p_i.
$$
Here $q_i$ denotes the original $Q$-component in $S$, which may be
nonzero even for $i\in A$.

\smallskip
\smallskip
\noindent
\hypertarget{main:selected-index-claim}{\textbf{A1.}} \emph{$F'_Q(x)=F_Q(x)$ for every
$x\in V$.}

\begin{proof}
Let $x\in V$. By the definitions of $F_Q$ and $F'_Q$, we have
$$
 F_Q(x)=F'_Q(x)+\sum_{\substack{i\in A\\ \mu_i(x)=1}}q_i.
$$
It therefore suffices to show that
$\sum\limits_{\substack{i\in A\\ \mu_i(x)=1}}q_i=0_Q$.
Summing \eqref{main:individual-Q-coordinate} over all $i\in A$ gives
the first equality below, since its right-hand side is zero when
$\mu_i(x)=0$. Interchanging the two finite sums then gives
$$
 \sum_{\substack{i\in A\\ \mu_i(x)=1}}q_i
  =\sum_{i\in A}
       \left(\sum_{\substack{\emptyset\ne I\subseteq J(x)\\ |I|\le a}}c_{i,I}\right)
  =\sum_{\substack{\emptyset\ne I\subseteq J(x)\\ |I|\le a}}
       \left(\sum_{i\in A}c_{i,I}\right)
 =0_Q \,,
$$
where the last equality follows from \eqref{main:zero-sum-coordinates}, whence $F'_Q(x)=F_Q(x)$. \qedhere \ ({\bf A1})
\end{proof}

For every $x\in V$, we have $F'_P(x)=F_P(x)$, since each index $i$ with
$\mu_i(x)=1$ contributes $p_i$ in both sequences. Together with {\bf A1}, this gives
$$
 F'(x)=F'_Q(x)+F'_P(x)=F_Q(x)+F_P(x)=F(x)
 \qquad\text{for every }x\in V,
$$
and hence
$$
 (F')^{-1}(0)=F^{-1}(0)=\{0,\mathbf1\},
$$
where the last equality follows from \eqref{main:zero-set}.
Thus \eqref{prelim:atom-criterion} shows that $S' \in \mathcal A (G_r)$. By construction, $|S'|=\ell$. The terms of $S'$ with
nonzero $Q$-component are precisely those indexed by
$i\in[1,\ell]\setminus A$ with $q_i\ne0_Q$. The number of such terms
is therefore at most $\ell-|A|=|CT^{-1}|$.
Since $CT^{-1} \in \mathcal F (Q^M)$ is zero-sumfree, we have
\begin{equation}\label{main:remaining-Q-bound}
|CT^{-1}|\le\mathsf d(Q^M)=M\mathsf d(Q) \,,
\end{equation}
where the equality follows from  \eqref{basic-1} and the subsequent comments.

For each $i\in[1,\ell]$, denote the $Q$-component of its $i$th term by $q_i'$, whence
$$
 q_i'=
 \begin{cases}
  0_Q & \text{if }i\in A,\\
  q_i & \text{if }i\in[1,\ell]\setminus A.
 \end{cases}    \qquad \text{and} \qquad  S'=\prod_{i=1}^\ell(g_i+q_i'+p_i).
$$
Put
$$
 I_0=\{i\in[1,\ell]\colon q_i'\ne0_Q\text{ or }g_i=0\}.
$$
For $k\in\N_0$, put
\begin{equation}\label{main:B-definition}
 B(k)=k+1+\mathsf d(Q)\sum_{s=1}^a\binom{k+1}{s}.
\end{equation}

\smallskip
\noindent
\hypertarget{main:index-bound-claim}{\textbf{A2.}} \emph{We have $|I_0|\le B(k_r)$.}

\begin{proof}
Let $e_1,\ldots,e_\ell$ be the standard basis of $\F_2^\ell$.
For every $i\in[1,\ell]$, we have $\varphi(e_i)=g_i$ by
\eqref{prelim:phi-definition}. Thus each $i$ with $g_i=0$ gives an element
$e_i\in\ker (\varphi)=V$. These elements are linearly independent, so
the dimension formula in \eqref{main:zero-set} gives
\begin{equation}\label{main:zero-binary-count}
 \bigl|\{i\in[1,\ell]\colon g_i=0\}\bigr|\le\dim V=\ell-r=k_r+1.
\end{equation}
Since $q_i'=0_Q$ for every $i\in A$, at most $\ell-|A|=|CT^{-1}|$
indices satisfy $q_i'\ne0_Q$. Combining this observation with
\eqref{main:zero-binary-count} and \eqref{main:remaining-Q-bound}, we obtain
$$
\begin{aligned}
 |I_0|&=\bigl|\{i\in[1,\ell]\colon q_i'\ne0_Q\text{ or }g_i=0\}\bigr|\\
 &\le \bigl|\{i\in[1,\ell]\colon q_i'\ne0_Q\}\bigr|
       +\bigl|\{i\in[1,\ell]\colon g_i=0\}\bigr|\\
 &\le |CT^{-1}|+k_r+1\\
 &\le M\mathsf d(Q)+k_r+1=B(k_r),
\end{aligned}
$$
where the last equality follows from $M=\sum_{s=1}^a\binom{k_r+1}{s}$. \qedhere \ ({\bf A2})
\end{proof}

\smallskip
We  construct a new group and a sequence $S''$ by replacing each term of $S'$ indexed by $I_0$
by two terms, using a separate new summand $C_2$ for each $i\in I_0$.
For each $i\in I_0$, adjoin a new direct summand
$\langle e_i'\rangle\cong C_2$. Thus the  enlarged group is
$$
 G_r\oplus\bigoplus_{i\in I_0}\langle e_i'\rangle
 \cong G_{r+|I_0|} \,,
$$
and we define
$$
 S''=\prod_{i\in[1,\ell]\setminus I_0}(g_i+q_i'+p_i)
      \prod_{i\in I_0}
      \bigl((g_i+e_i'+p_i) (e_i'+q_i')\bigr).
$$

\smallskip
\noindent\hypertarget{main:double-prime-minimality}{\textbf{A3.}}
\emph{$S''$ is a minimal zero-sum sequence over $G_{r+|I_0|}$ of length $\ell+|I_0|$.}

\begin{proof}
Each of the $|I_0|$ replacements adds one term, so $|S''|=\ell+|I_0|$.
For every $i\in I_0$, the two new terms have sum
$(g_i+e_i'+p_i)+(e_i'+q_i')=g_i+q_i'+p_i.$
Consequently, $\sigma(S'')=\sigma(S')=0$.
Consider a nontrivial zero-sum subsequence $S^*$ of $S''$. Then the
projection onto $\langle e_i'\rangle$ shows that either both terms of the
$i$th pair are selected or neither, for each $i\in I_0$.
Contracting each selected pair to $g_i+q_i'+p_i$ therefore turns $S^*$ into a nontrivial zero-sum subsequence
of $S'$. By the minimality of $S'$, the contracted sequence equals $S'$.
Hence both new terms are selected for every $i\in I_0$, and the original
term is selected for every $i\in[1,\ell]\setminus I_0$.
Thus  $S^*=S''$. \qedhere \ ({\bf A3})
\end{proof}

\smallskip
\noindent\textbf{A4.} \emph{For every integer $j\ge |I_0|$, there exists a
minimal zero-sum sequence over $G_{r+j}$ of length $\ell+j$ such that the order
of each of its terms divides either $n/2$ or $2^a$.}

\begin{proof}
First we prove the claim for  $j=|I_0|$,  and we assert that  $S''$ has the required properties.
By \hyperlink{main:double-prime-minimality}{\textbf{A3.}}, it remains to check the orders of its terms.
For $i\in[1,\ell]\setminus I_0$, the term in the first product
defining $S''$ is $g_i+q_i'+p_i=g_i+p_i$, since $q_i'=0_Q$.
Here $g_i\ne0$ has order $2$, and $\ord(p_i)$ divides the odd integer $p^b$.
Since $g_i$ and $p_i$ lie in distinct summands, we obtain
$$
 \ord(g_i+p_i)=2\ord(p_i)\mid 2p^b.
$$

For $i\in I_0$, consider the two terms of the $i$th pair in the second product.
For the first term, $\ord(g_i+e_i')=2$ because $g_i+e_i'$
lies in the elementary $2$-group and has nonzero component $e_i'$ in the new
summand. Since $p_i$ lies in $P$ and $\ord(p_i)$ divides the odd integer $p^b$,
we have
$$
 \ord(g_i+e_i'+p_i)
   =\lcm\bigl(\ord(g_i+e_i'),\ord(p_i)\bigr)
   =\lcm\bigl(2,\ord(p_i)\bigr)
    =2\ord(p_i)\mid 2p^b.
$$
For the second term, $e_i'$ and $q_i'$ lie in distinct summands.
Since $\ord(e_i')=2$ and $\ord(q_i')\mid 2^a$, we obtain
$$
\ord(e_i'+q_i')=\lcm\bigl(2,\ord(q_i')\bigr)\mid 2^a.
$$
Since $a\ge2$, we have $2p^b\mid 2^{a-1}p^b=n/2$.
This proves the case $j=|I_0|$.

Now we consider the case $j>|I_0|$. Put $t=j-|I_0|$ and choose a term $h$ of
$S''$ with $\ord(h)\mid n/2$. If $I_0\ne\emptyset$, take $h=g_i+e_i'+p_i$
for some $i\in I_0$. If $I_0=\emptyset$, take $h$ to be any term $g_i+p_i$ of $S''$.

Adjoin $t$ new direct summands $\langle w_1\rangle,\ldots,\langle w_t\rangle$,
each of order $2$. Identify
$G_{r+|I_0|}\oplus\bigoplus_{i=1}^t\langle w_i\rangle$ with $G_{r+j}$,
and regard $S''$ as a sequence over the first summand. Define
$$
 \widetilde{S''}
   =(S''h^{-1})\cdot(h+w_1+\ldots+w_t)\cdot w_1\ldots w_t
   \ \in\mathcal F(G_{r+j}).
$$ Since $\sigma(S'')=0$, we can verify that $\widetilde{S''}$ is a minimal zero-sum sequence of length $|\widetilde{S''}|=|S''|-1+(t+1)=\ell+j.$

Finally, each $w_i$ also has order $2$, and so
$$
 \ord(h+w_1+\ldots+w_t)=\lcm\bigl(\ord(h),2\bigr)\mid n/2.
 \eqno\hbox{\qedhere \ ({\bf A4})}
$$
\end{proof}

\medskip

For $j\ge |I_0|$, {\bf A4} gives a minimal zero-sum sequence over
$G_{r+j}$ with the required term orders and length $\ell+j=\mathsf D(G_r)+j$.
Theorem~\ref{inv:main} requires length $\mathsf D(G_{r+j})$. By \eqref{main:k-definition}, we have
$$
 \mathsf D(G_{r+j})-(\ell+j)=k_{r+j}-k_r.
$$
Thus the minimal zero-sum sequence from {\bf A4} has maximal length equal to $\mathsf D(G_{r+j})$ precisely when $k_{r+j}=k_r$.

\medskip
Now, we call a rank $s\in\N_0$ \emph{exceptional} if every $S \in \mathcal A (G_s)$ of length $\mathsf D(G_s)$
has a term $g$ such that
$$
 \ord(g)\nmid n/2 \qquad\text{and}\qquad \ord(g)\nmid 2^a.
$$
\smallskip
\noindent\hypertarget{main:exceptional-ranks-claim}{\textbf{A5.}}
\emph{For every $k\in\N_0$, at most $B(k)$ (see \eqref{main:B-definition}) exceptional ranks $s$ satisfy $k_s=k$.}

\begin{proof}
Fix $k\in\N_0$. If no rank $s$ satisfies $k_s=k$, the count is zero. Otherwise put
$$
 r_0=\min\{s\in\N_0\colon k_s=k\}.
$$
By \eqref{main:rank-monotonicity}, $k_{r_0+B(k)}\ge k_{r_0}=k$.
If $k_{r_0+B(k)}>k$, then $k_s>k$ for every $s\ge r_0+B(k)$.
Hence every rank $s$ with $k_s=k$ lies in $[r_0,r_0+B(k)-1]$,
which contains exactly $B(k)$ ranks. We may therefore assume that
$k_{r_0+B(k)}=k$.

Take $r=r_0$ and choose a minimal zero-sum sequence $S$ over
$G_{r_0}$ of length $\ell=\mathsf D(G_{r_0})$, as at the beginning of the
proof. Construct $S'$ from $S$ as in \eqref{main:modified-sequence}, and let $I_0$
be the set of indices of the terms of $S'$ with $q_i'\ne0_Q$ or $g_i=0$. Since $k_{r_0}=k$, \hyperlink{main:index-bound-claim}{A2} gives $$|I_0|\le B(k_{r_0})=B(k).$$

Let $s\ge r_0+B(k)$ be any rank with $k_s=k$, and put
$j=s-r_0$. Then $j\ge B(k)\ge |I_0|$.
{\bf A4} applied with $r=r_0$ therefore gives a minimal zero-sum
sequence over $G_{r_0+j}=G_s$ of length $\ell+j$ satisfying the required
order conditions.
By \eqref{main:k-definition} with $r=s$, we have
$\mathsf D(G_s)=s+k_s+1$.
Since $k_s=k_{r_0}=k$, the constructed sequence has length
$$
 \ell+j
 =\mathsf D(G_{r_0})+(s-r_0)
 =(r_0+k+1)+(s-r_0)
 =s+k+1
 =\mathsf D(G_s).
$$
Thus any exceptional rank $s$ with $k_s=k$ must satisfy
$r_0\le s<r_0+B(k)$. There are at most $B(k)$ such ranks. \qedhere \ ({\bf A5})
\end{proof}

For $x\ge0$, put $K=k_{\lfloor x\rfloor}$. Since
$(k_s)_{s\ge0}$ is nondecreasing by \eqref{main:rank-monotonicity}, every rank $s\in[0,\lfloor x\rfloor]$
satisfies $0\le k_s\le K$. Group the exceptional ranks $s\in[0,\lfloor x\rfloor]$ by
the value of $k_s$. Each such rank $s$ is counted exactly once
on the right-hand side of \eqref{main:exceptional-count}, in the summand indexed by $k=k_s$.
\begin{equation}\label{main:exceptional-count}
 \bigl|\{s\in[0,\lfloor x\rfloor]\colon s\text{ is exceptional}\}\bigr|
 =\sum_{k=0}^K
 \bigl|\{s\in[0,\lfloor x\rfloor]\colon k_s=k,\ s\text{ is exceptional}\}\bigr|.
\end{equation}
By \hyperlink{main:exceptional-ranks-claim}{A5}, the summand for each $k$ is at most $B(k)$.

Fix $k\in\N_0$. For $1\le t\le a$, we have
$\binom{k+1}{t}\le(k+1)^t\le(k+1)^a$.
Applying this bound to the $a$ summands in \eqref{main:B-definition} gives
$$
 B(k)
 \le k+1+\mathsf d(Q)\sum_{t=1}^a(k+1)^a
 = k+1+a\mathsf d(Q)(k+1)^a
  \le (1+a\mathsf d(Q))(k+1)^a,
$$
where the last inequality uses $k+1\le(k+1)^a$, since $k+1\ge1$ and $a\ge2$.
Hence we obtain that
$$
 \bigl|\{s\in[0,\lfloor x\rfloor]\colon s\text{ is exceptional}\}\bigr|
 \le \sum_{k=0}^K B(k)
 \le (1+a\mathsf d(Q))\sum_{k=0}^K(k+1)^a
 \le (1+a\mathsf d(Q))(K+1)^{a+1}.
$$
Applying \eqref{main:k-logarithmic} with $r=\lfloor x\rfloor$
gives $K=k_{\lfloor x\rfloor}\ge n-1\ge1$ and
$K\ll_{G_0}\log(\lfloor x\rfloor+2)$. Hence
$K+1\le 2K\ll_{G_0}\log(x+2)$. Thus the number of exceptional ranks is
$O_{G_0}((\log(x+2))^{a+1})$. Since $a$ is fixed,
$(\log(x+2))^{a+1}/(\lfloor x\rfloor+1)\to0$ as $x\to\infty$.
The exceptional ranks therefore have natural density zero, and the set
of ranks in Theorem~\ref{inv:main} has natural density one.
\end{proof}

\section{Proof of Theorem~\ref{inv:general-counterpart}}\label{sec:general-counterpart}

Let $G$ be a finite abelian group. Let $\eta(G)$ be the smallest integer $\ell$ such that every sequence $S \in \mathcal F (G)$ of length $|S| \ge \ell$ has a  zero-sum subsequence $T$ of length  $|T| \in [1, \exp(G)]$.

\smallskip
The $\eta (G)$-invariant is a classic zero-sum invariant studied since the 1960s. Among others, we have
\begin{itemize}
\item[(i)] If $G = C_{n_1} \oplus C_{n_2}$ with $1 \le n_1 \mid n_2$, then $\eta (G) = 2n_1+n_2-2$ (\cite[Theorem 5.8.3]{Ge-HK06a}).

\item[(ii)] If $G$ is a $p$-group with $\mathsf D (G) \le 2 \exp (G)-1$, then
      \begin{equation} \label{eta-precise}
      \eta (G) = 2 \mathsf D (G) - \exp (G) \qquad \text{(\cite[Theorem~1.6]{Lu17})} \,.
      \end{equation}
\end{itemize}
The precise value of $\eta (C_p \oplus C_p \oplus C_p)$ is known only for  some small primes $p \in \mathbb P$.

\medskip
\begin{proof}[Proof of Theorem~\ref{inv:general-counterpart}]
We have $G_k=G\oplus C_{kp^m}\cong C_{n_1}\oplus\ldots\oplus C_{n_r}$ with $1<n_1\mid\ldots\mid n_r$.
If $|G|=1$, then $G_k$ is cyclic and both assertions hold. Suppose that $|G|>1$, furthermore,  we may suppose without restriction that $p \nmid k$. Since $\exp (G) \le \mathsf D (G) \le p^m$, it follows that $\exp (G_k)=kp^m$.  By \cite[Theorem 3.3.10]{Ge-Gr-Zh26a}, we have $\mathsf D (G_k) = kp^m + \mathsf d (G)$.
Let $\varphi:G_k = G \oplus C_{kp^m} \to H$ be the identity on $G$ and the multiplication by $k$ on the second summand. Thus, $H \cong G\oplus C_{p^m}$ is a $p$-group, $\ker (\varphi)$ is cyclic of order $k$,  $\mathsf D(H)=p^m+\mathsf d (G)$, and $\eta(H)=p^m+2\mathsf d (G)$ by \eqref{eta-precise}. Note, if we have $H \cong C_{m_1} \oplus \ldots \oplus C_{m_{r'}}$ with $1 < m_1 \mid \ldots \mid m_{r'}$, then $m_1=n_1$.
Let  $S=g_1 \cdot\ldots \cdot g_{kp^m+\mathsf d (G)} \in \mathcal A (G_k)$ of length $|S|= \mathsf D (G_k)$.

\medskip
\noindent
{\bf Step 1.}
We consider  decompositions of the form
\begin{equation} \label{basic-decomposition}
S=U_1\ \cdot \ldots \cdot U_s  R \,,
\end{equation}
where $s \in [0,k]$ ($R=S$ when $s=0$), $U_1, \ldots, U_s, R \in \mathcal F (G_k)$ with $\varphi(U_i)\in\mathcal A(H)$ for each $i \in [1,s]$, and we set $L=\sum_{i=1}^s|U_i|=|S|-|R|$. All parameters $s, R, L$ and the forthcoming integer $b$ depend on the given decomposition. In this first step, we derive some properties of these parameters, valid for any such decomposition. In the second step we use this for different decompositions of the form \eqref{basic-decomposition}.

Since $\sigma(\varphi(R))=\varphi(\sigma(S))-\sum_{i=1}^s\sigma(\varphi(U_i))=0$, $R$ allows a decomposition $R=V_1 \cdot\ldots \cdot V_b$, where $V_1, \ldots, V_b \in \mathcal F (G_k)$ such that  $\varphi(V_j)\in\mathcal A(H)$ for  $j \in [1,b]$, and suppose that $b$ maximal such that a decomposition with these properties exists. Then
$$
 \sigma(U_1) \cdot\ldots \cdot\sigma(U_s) \sigma(V_1) \cdot\ldots \cdot\sigma(V_b) \ \in \mathcal F \big( \ker (\varphi) \big)
$$
is a  zero-sum sequence. Moreover, it is minimal because any proper nontrivial zero-sum subsequence of it would give a proper nontrivial zero-sum subsequence of $S$. Thus $s+b\le \mathsf D ( \ker (\varphi)) = \mathsf D (C_k) = k$.

\smallskip
We assert that  $|R|\le bp^m+\mathsf d (G)$. If $b=0$, then $R=1$. Otherwise, assume to the contrary that $b \ge 1$ and that  $|R|>bp^m+\mathsf d (G)$. Choose some element $f\in \supp ( R)$. Since $\mathsf D(H)=p^m+\mathsf d (G)<2p^m$, \cite[Theorem 6.1.5]{Ge-HK06a} gives $b$ disjoint nontrivial zero-sum subsequences of $\varphi(R  f^{-1})$. Thus we obtain $W_1 \cdot\ldots \cdot W_b\mid R  f^{-1}$ with $\varphi(W_j)\in\mathcal A(H)$. Then $R (W_1 \cdot\ldots \cdot W_b)^{-1}$ contains $f$ and has zero-sum image, so factoring its image gives more than $b$ factors of $R$, a contradiction.

Therefore, it follows that
$$
 kp^m+\mathsf d (G)-L=|R|\le bp^m+\mathsf d (G) \qquad\text{and}\qquad k-L/p^m\le b\le k-s.
$$
In particular, $L\ge sp^m$, and $L<(s+1)p^m$ forces $s+b=k$.

\medskip
\noindent
{\bf Step 2.}
We choose some  $g\in \supp ( S)$ and will find a decomposition \eqref{basic-decomposition} having additional properties.
Let $i\in[1,k-1]$ and suppose that $U_1,\ldots,U_{i-1}$ have already been chosen with lengths at most $p^m$. Then the sequence $Sg^{-1}(U_1\cdot \ldots \cdot U_{i-1})^{-1}$ has length at least $2p^m+\mathsf d(G)-1\ge\eta(H)$. Thus we can choose a subsequence $U_i$ of this remainder with $\varphi(U_i)\in\mathcal A(H)$ and $|U_i|\le p^m$.  For $s=k-1$, the bound $L\ge(k-1)p^m$ forces $|U_i|=p^m$ for all $i$. Put
\begin{equation}\label{decomposition-with-T}
 T=S(U_1\cdot \ldots \cdot U_{k-1})^{-1}.
\end{equation}
Then $g \mid T$  and $|T|=p^m+\mathsf d(G)$. Since $L=(k-1)p^m<kp^m$, we have $s+b=k$, whence $b=1$ and $\varphi(T)\in\mathcal A(H)$.

Next we choose some  $f \in \supp( T  g^{-1})$. Since $\varphi (T) \in \mathcal A (H)$, the sequence $\varphi(T  f^{-1})$ is zero-sum free of length $\mathsf d(H)$.  By \cite[Theorem~5.3]{Gi10a}, we infer that $m_1\mid\ord(\varphi(g))$.  Since $m_1=n_1$ and $\ord ( \varphi (g) ) \mid \ord (g)$, the assertion of  Part~2 follows.   It remains to prove Part~1.

\smallskip
By \cite[Corollary~2.3]{Gi10a}, $\varphi(T f^{-1})$ contains at least $p^m-1$ terms of order $p^m$. Choose some $f' \in \supp (T f^{-1})$ with $\ord \big(\varphi(f') \big)=p^m$. The sequence $\varphi(T {f'}^{-1})$ is zero-sum free of length $\mathsf d(H)$. Again \cite[Corollary~2.3]{Gi10a} shows that there are at least $p^m-1$ terms of order $p^m$ in this sequence. Thus $\varphi (T)$ has  at least $p^m$ terms of order $p^m$.
If $k=1$, then $\varphi$ is an isomorphism, giving at least $p^m\ge p^m-\mathsf d (G)$ terms of order $kp^m=p^m$, as required.
If $k > 1$, then we distinguish two cases.

\smallskip
\noindent
CASE 1:  $k$ is a power of a prime $q$.

We set $K=G\oplus q C_{kp^m}$ and consider the canonical epimorphism $\psi:G_k\to G_k/K$. Let $u \in \N_0$ be the number of terms of $S$ in $K$. Then we choose a decomposition $S=Z_1\cdot\ldots\cdot Z_v$, where $|Z_1|= \ldots =|Z_u|=1$ and   $\psi(Z_i)\in\mathcal A(G_k/K)$ for all $i \in [u+1, v]$. Note that $|Z_i| \le q$  for all $i \in [u+1,v]$. Then $\sigma(Z_1) \cdot\ldots \cdot\sigma(Z_v)\in\mathcal A(K)$,  $\mathsf D(K)=kp^m/q+\mathsf d(G)$ by \cite[Theorem 3.3.10]{Ge-Gr-Zh26a}, and hence
$$
 u+(kp^m+\mathsf d (G)-u)/q = u + (|S|-u)/q \le v\le  \mathsf D(K)=kp^m/q+\mathsf d(G)\,.
$$
This gives $(1-1/q)u\le(1-1/q)\mathsf d(G)$, so $u\le\mathsf d(G)$. Of the at least $p^m$ terms $f$ of $T$ from \eqref{decomposition-with-T} with $\ord(\varphi(f))=p^m$, at most $u$ lie in $K$, since $T\mid S$. Thus at least $p^m-\mathsf d(G)$ lie outside $K$.
Let $f$ be such an element, say $f = h + ae \in G \oplus \langle e \rangle = G \oplus C_{kp^m}$ with
 $q\nmid a$. Hence $k\mid\ord(ae)\mid\ord(f)$ and $p^m=\ord(\varphi(f))\mid\ord(f)$. Since $\gcd(k,p^m)=1$ and $\ord(f)\mid kp^m$, we obtain $\ord(f)=kp^m$.

\medskip
\noindent
CASE 2: $k$ is divisible by two distinct primes.

Note that  $k\ge6$.  We start with the following assertion.

\smallskip
\noindent
{\bf A1.} If $U$ and $V$ are subsequences of $S$ such that $\varphi (U), \varphi (V) \in \mathcal A (H)$, then $\sigma (U)=\sigma (V)$ and $\ord \big( \sigma (U) \big)  = k$. We set $c:= \sigma (U)$.

\noindent
\begin{proof}[Proof of {\bf A1}]
Let $U$ and $V$ be subsequences of $S$ with $\varphi(U),\varphi(V)\in\mathcal A(H)$, say $U = \prod_{i\in I} g_i$ and $V= \prod_{j \in J}g_j$ and set  $R_{U,V}=\prod_{i\in[1,kp^m+\mathsf d (G)]\setminus(I\cup J)}g_i$. Then  $|R_{U,V}| = |S|-|U|-|V|+|I\cap J| \ge |S|-|U|-|V|$. Note that  $|U| =|\varphi (U)| \le \mathsf D(H)=p^m+\mathsf d(G)$ and, similarly, $|V| \le \mathsf D(H)=p^m+\mathsf d(G)$. Since $k\ge6$ and $\mathsf d (G)<p^m$, we have
$$
 |R_{U,V}|\ge |S|-2(p^m+\mathsf d (G))=(k-2)p^m-\mathsf d (G)\ge4p^m-\mathsf d (G)>p^m+2\mathsf d (G)=\eta(H).
$$
By the definition of $\eta(H)$,  $R_{U,V}$ has a subsequence $W$ with $\varphi(W)\in\mathcal A(H)$ and $|W|\le p^m$. We apply \eqref{basic-decomposition}, with $s=1$, $U_1=W$, and $R=SW^{-1}$. Then we obtain that $|W|=L\ge p^m$, whence $|W|=p^m$.

Next apply \eqref{basic-decomposition} with $s=2$ to $(U,W)$ and $(V,W)$ separately, with remainders $S(UW)^{-1}$ and $S(VW)^{-1}$. In each case $L\le\mathsf D(H)+p^m<3p^m$, so both pairs extend to factorizations into $k$ factors, whose $\varphi$-images are minimal zero-sum sequences.

In each factorization the factor sums form a length-$k$ minimal zero-sum sequence over $\ker(\varphi)\cong C_k$ and hence all equal a generator. The common factor $W$ gives $\sigma(U)=\sigma(W)=\sigma(V)$.
Thus all subsequences $U$ of $S$ with $\varphi (U)\in\mathcal A(H)$ have the same sum, say $\sigma (U) = c$ with  $\ord (c)=k$. A decomposition of $S$ into $\ell$ such subsequences gives $\ell c=\sigma(S)=0$ and $\ell \in [1, k]$, whence $\ell=k$.
\qedhere \ ({\bf A1})
\end{proof}

\smallskip
We  set $A=ST^{-1}$ with $T$ from \eqref{decomposition-with-T}. Then $|A|=(k-1)p^m$ and $\varphi(A)$ is zero-sum. If $U$ is a subsequence of $A$ such that $\varphi (U)$ has sum zero, then  $\varphi (AU^{-1})$ also has sum zero. Decompose both into factors whose $\varphi$-images are minimal zero-sum sequences. Adjoining $T$ gives $k$ factors in all. The $k-1$ factors outside $T$ each have length at least $p^m$ (take $s=1$) and total length $|A|=(k-1)p^m$, so each has length $p^m$ and sum $c$. Hence
$$
 p^m\mid |U| \qquad\text{and}\qquad \sigma(U)=(|U|/p^m)c.
$$

\medskip
 We use the group algebra $\F_p[H]$ and all notation, introduced in \eqref{eq:augmentation-ideal}, \eqref{vanishing-of-the-augmentation-ideal}, and \eqref{basic-notation}.
Choose a decomposition $A=B D R$, with $B,D, R \in \mathcal F (G_k)$, that  maximizes $|B|+|D|$ subject to the condition that
$$
\alpha(y)\beta(-y)\ne 0 \ \text{ for some $y\in H$, where} \  \alpha(y)=[X^y]\Pi \big(\varphi (B) \big) \ \text{ and } \ \beta(y)=[X^y]\Pi \big(\varphi (D) \big) \,.
$$
Such a decomposition exists, since $A$ has finitely many decompositions only and taking $B=D=1_{\mathcal F (G_k)}$  gives $\alpha(0)\beta(0)=1$. We continue with the following assertion.

\smallskip
\noindent
{\bf A2.} Every $g \in \supp (R)$ has order $kp^m$ and $\max \{|B|, |D|\} \le \mathsf d (H)$.

Suppose that {\bf A2} holds and recall that $\mathsf d (H) = p^m + \mathsf d (G)-1$. Then we have $|R|$ elements of order $kp^m$ and
$$
\begin{aligned}
|R| & = |A| -|B|-|D| = (k-1)p^m - |B| - |D| \\
 & \ge (k-1)p^m - 2(p^m + \mathsf d (G) -1) \ge p^m - \mathsf d (G) \,.
\end{aligned}
$$
Thus, it remains to prove {\bf A2.}

\noindent
\begin{proof}[Proof of {\bf A2}]
We choose $g \in \supp ( R)$ and set $x=\varphi(g)$. The factorizations $A=(B  g) D (R g^{-1})$ and $A=B (D  g) (R  g^{-1})$ violate the coefficient condition by maximality. Since $\Pi \big(\varphi (B g) \big)=(1-X^x)\Pi \big(\varphi (B) \big)$ and $\Pi \big(\varphi (D  g) \big)=(1-X^x)\Pi \big(\varphi (D) \big)$, taking the coefficient pairs at $y$ and $y-x$, respectively, gives for all $y\in H$
$$
 \alpha(y)\beta(-y)=\alpha(y-x)\beta(-y)=\alpha(y-x)\beta(-y+x).
$$

Choose $y_0$ with $\alpha(y_0)\beta(-y_0)\ne0$ and put $y_j=y_0+jx$ for $j\in[0,t-1]$, where $t=\ord(x)$. The preceding identities make $\alpha(y_j)\beta(-y_j)\ne0$ for every $j \in [0,t-1]$. Expanding $\Pi \big(\varphi (B) \big)$ and $\Pi \big(\varphi (D) \big)$ gives $B_j\mid B$ and $D_j\mid D$ with $\sigma \big(\varphi(B_j) \big)=y_j$ and $\sigma \big(\varphi(D_j) \big)=-y_j$. With $B_{-1}=B_{t-1}$, both $B_j D_j$ and $B_{j-1} g D_j$ are subsequences of $A$ such that $\varphi (B_jD_j)$ and $\varphi (B_{j-1}gD_j)$ have sum zero for every $j\in[0,t-1]$. Hence these subsequences satisfy
$$
 p^m\mid |B_j  D_j|,\quad p^m\mid |B_{j-1}  g  D_j|,\quad
 \sigma(B_j  D_j)=\frac{|B_j  D_j|}{p^m}c\quad\text{and}\quad
 \sigma(B_{j-1} g D_j)=\frac{|B_{j-1} g D_j|}{p^m}c.
$$
Setting
$$
 E=\prod_{j=0}^{t-1}(B_j D_j)\qquad\text{and}
 \qquad E'=\prod_{j=0}^{t-1}(B_{j-1} g D_j)=g^t  E
$$
we obtain that
$$
 p^m\mid |E|,\quad p^m\mid |E'|,\qquad
 \sigma(E)=\frac{|E|}{p^m}c\quad\text{and}\quad\sigma(E')=\frac{|E'|}{p^m}c.
$$
Since $|E'|-|E|=t$, we infer that $p^m\mid t$, and thus $t=p^m$ because $\exp(H)=p^m$. Now $p^mg=\sigma(E')-\sigma(E)=c$. Together with $\ord \big(\varphi(g) \big)=p^m$ and $\ord(c)=k$, this yields $\ord(g)=kp^m$.

Both $\Pi \big(\varphi (B) \big)$ and $\Pi \big(\varphi (D) \big)$ are nonzero, and each factor in either product lies in the augmentation ideal $I_H \subseteq \F_p[H]$. Since  $I_H^{\mathsf d(H)+1}=0$ by \eqref{vanishing-of-the-augmentation-ideal}, it follows that  $|B| \le\mathsf d(H)$ and $|D|\le\mathsf d(H)$. \qedhere \ ({\bf A2})
\end{proof}

\end{proof}

\section{Proof of Theorem~\ref{inv:girard-result}}\label{sec:girard}

The following lemma will be used as a technical tool in the proof of Theorem~\ref{inv:girard-result}.

\begin{lemma}[\cite{Ga-Li-Zh-Zh16}, Lemma~4.2]\label{girard:binary-blocks}
	Let $\ell,m, u \in \N$  with $m\ge2$ and $\ell \in [2, m]$ such that
	$$
	2^{u-1}\mid\ell
	\qquad\text{and}\qquad
	\ell\frac{2^u-1}{2^{u-1}}-u\le m.
	$$
	Then there exists a sequence $S \in \mathcal F (C_2^m)$  of length $|S|=m+u$ whose nontrivial zero-sum subsequences all have length $\ell$.
\end{lemma}

\begingroup
\newcommand{\girardcodeD}{\mathord{\text{\usefont{U}{rsfs}{m}{n}D}}}


For an element $x = (x_1+2\Z, \ldots, x_r+2\Z) \in \F_2^r$, with $x_1, \ldots, x_r \in \{0,1\}$, the {\it Hamming weight} $\mathsf{wt} (x)$ is defined as $\mathsf{wt}(x) = \sum_{i=1}^r x_i \in [0,r]$.
By the dimension of an affine subspace we mean the dimension of the associated linear space.  For positive integers $r$ and $M$, let $\girardcodeD(r,M)$ be the largest dimension of an affine subspace $C\subseteq\F_2^r$ for which there exists some $a \in [0,M]$ such that
$$
\bigl|\{x\in C\colon \mathsf{wt}(x)\equiv a \mod M\}\bigr|=1
$$
(this is a concept introduced in \cite[Definition~4.3]{Br-Go-Gu22}).
Next we recall three results (partly in a simpler form which is sufficient for our needs) that will be used in the proof of Part~1 of Theorem~\ref{inv:girard-result}.

\begin{lemma}[\cite{Go-Gr-Ma-Ta26}, Theorem~1.1]\label{girard:torsion-pfr}
	Let $G$ be a  finite  abelian group and let  $A\subseteq G$ be a finite nonempty subset. If $|A+A|\le K|A|$ for some real number $K >  0$, then there is a subgroup $H\le G$ with  $|H| \le |A|$ such that at most $(2K)^{O(\exp (G)^3)}$ cosets of $H$ cover $A$.
\end{lemma}

\begin{lemma}[Pl\"unnecke's Inequality]\label{girard:plunnecke}
Let $G$ be an abelian group and let $A \subseteq G$ be a finite nonempty subset. If $|A-A| \le K|A|$ for some real number $K > 0$, then $|A+A| \le K^2|A|$.
\end{lemma}

\begin{proof}
This follows from  \cite[Corollary 2.2.3]{Ge-Gr-Zh26a} (with $B = -A$).
\end{proof}

%

``Assuming the {\it Polynomial Freiman--Ruzsa (PFR) Conjecture}'' means the following (see \cite{Bh-Dv-Lo14}). For every fixed integer $q\ge2$, there must exist a constant $\gamma_q>0$, depending only on $q$, such that for every $t\in\N$, every real $K\ge2$, and every nonempty subset $A\subseteq Q_t=(\Z/q\Z)^t$ satisfying $|A-A|\le K|A|$, there is a subgroup $H\le Q_t$ with $|H|\le |A|$ such that $A$ is covered by at most $K^{\gamma_q}$ cosets of $H$. In particular, $\gamma_q$ is independent of $t$, $K$, and $A$.


\begin{lemma}[\cite{Br-Go-Gu22}, Corollary~4.24]\label{girard:code-bound}
	Assuming the PFR conjecture in the above sense, for every integer $q\ge2$, there exists a constant $c_q>0$ depending only on $q$ such that, for every integer $n\ge2$,
	$$
	\girardcodeD(n,q)\leqslant c_q\frac{n}{\log n}.
	$$
\end{lemma}

We cite the  corollary in its original form, but (in our setting) the hypothesis is now a consequence of Lemma~\ref{girard:plunnecke} and Lemma~\ref{girard:torsion-pfr}. Indeed, fix any $q\ge2$ and let $t$, $K$, and $A$ be as above. Lemma~\ref{girard:plunnecke} yields $|A+A|\le K^2|A|$. Since $\exp(Q_t)=q$, Lemma~\ref{girard:torsion-pfr} gives a subgroup $H\le Q_t$ with $|H|\le |A|$ whose translates cover $A$, with the number of translates at most
$$
 (2K^2)^{c q^3}\le K^{3c q^3}\qquad(K\ge2),
$$
where $c>0$ is an absolute constant. Thus $\gamma_q=3c q^3$ is admissible, uniformly for every $t\in\N$. Since $q$ was arbitrary, the required PFR covering statement holds for every fixed modulus $q\ge2$, and Lemma~\ref{girard:code-bound} can be applied unconditionally.

\medskip
\begin{proof}[Proof of Theorem~\ref{inv:girard-result}]
	Let  $n \ge 3$ be odd and let $r \in \N$.  We think of $n$ being fixed and $r$ being a growing parameter.
	We set $G_r=\F_2^{r+1}\oplus \Z/n\Z$. If $h=g+a \in G_r$, with $g \in \F_2^{r+1}$ and $a \in \Z/n\Z$, then $\ord (h) = \lcm ( \ord (g), \ord (a) )$. In particular, $h$ has odd order if and only if $g=0$. We call $r$ {\it exceptional} if there exists a zero-sum free sequence of length $\mathsf d(G_r)$ having an element of odd order.
Equivalently, $G_r$ admits a minimal zero-sum sequence of length $\mathsf D(G_r)$ containing an element of odd order.

	By \eqref{basic-1}, we have $\mathsf d (G_r) \ge r+1$ and we set
	$k_r=\mathsf d(G_r)-r$. Since $\mathsf d (G_{r+1}) \ge \mathsf d (G_r) +1 $, we have
	\begin{equation} \label{k_r-inequality}
		k_{r+1} \ge k_r \,.
	\end{equation}
	Let
	\begin{equation} \label{given-sequence}
		S=\prod_{i=1}^{\mathsf d (G_r)} (g_i+a_i) \in \mathcal F (G_r)
	\end{equation}
	be a zero-sum free sequence, where $g_i\in\F_2^{r+1}$ and $a_i\in\Z/n\Z$.
	
	\smallskip
	\noindent\emph{The block construction.}
	For $j\ge1$, put
	$$
	r_j=2^j-1-j,\qquad \ell_j=2^{j-1},\qquad M_j=r_j+j=2^j-1.
	$$
	We construct a sequence $T_j \in \mathcal F \big( \F_2^{r_j} \big)$ of length $M_j$  whose nontrivial zero-sum subsequences all have length $\ell_j$. For $j\ge3$, apply Lemma~\ref{girard:binary-blocks} with $(m,\ell,u)=(r_j,\ell_j,j)$. To check the hypotheses, note that $2^{j-1}\ge j+1$ for $j\ge3$: equality holds at $j=3$, and doubling the left side at each step preserves the inequality. Thus
	$$
	r_j-\ell_j=2^{j-1}-j-1\ge0,\qquad
	\ell_j=2^{u-1}\ge2,\qquad
	\ell_j\frac{2^u-1}{2^{u-1}}-u=2^j-1-j=r_j.
	$$
	This gives $m\ge2$, $\ell \in [2, m]$, the required divisibility, and the rank bound. The lemma supplies a sequence of length $m+u=M_j$. For $j=1$, take the  zero  in $\F_2^0$, so $M_1=\ell_1=1$. For $j=2$, take three copies of the nonzero element $e$ in $\F_2^1$. A nontrivial subsequence sums to zero exactly when it selects two copies, since $e+e=0$ and the other possible lengths are $1$ and $3$. This gives $M_2=3$ and $\ell_2=2$.
	
	For $a\in \Z / n \Z$, let $b \in \Z/n\Z$ be the unique element with $\ell_j b=a$. Such an element exists and is unique because $\gcd(\ell_j,n)=1$. Writing $T_j=t_1\cdot \ldots \cdot  t_{M_j}$, we define
	$$
	B_j(a)=\prod_{h=1}^{M_j}(t_h + b)
	\ \in \ \mathcal F \big( \F_2^{r_j}\oplus \Z / n \Z \big).
	$$
	A nontrivial subsequence of $B_j(a)$, whose first component sums to zero,   selects a nontrivial zero-sum subsequence of $T_j$ at the same positions. It therefore has exactly $\ell_j$ terms. Each of them has second component $b$, whence their sum in $\Z / n \Z$ is $\ell_j b=a$.
	
	We next combine blocks. Let $s \in [1,n-1]$ and suppose $\sum_{i=1}^s r_{j_i}\le r+1$. Choose subspaces $U_1,\ldots,U_s,U_0$ of $\F_2^{r+1}$ such that
	$$
	\F_2^{r+1}=U_1\oplus\ldots\oplus U_s\oplus U_0,\qquad
	\dim U_i=r_{j_i}\ (i \in [1, s]),\qquad
	\dim U_0=r+1-\sum_{i=1}^s r_{j_i}.
	$$
	Such a choice is possible. For example, partition a basis of $\F_2^{r+1}$ into sets of these sizes and take their spans. Embed the first component space of $B_{j_i}(1+n\Z)$ as $U_i$. Form a sequence from all the blocks and one term $(w,0)$ for each basis vector $w$ of $U_0$.
	
	Consider a nontrivial subsequence of the sequence just constructed whose first component sums to zero. Its component in each direct summand must vanish. In $U_0$, only the added basis terms can contribute, and no nonempty subset of a basis sums to zero. Hence the subsequence selects no added term. In $U_i$, only the terms of $B_{j_i}(1)$ can contribute. If the subsequence selects any terms in that block, those terms have first-component sum zero and contribute $1+n\Z$ to the second-component sum. Thus the second-component sum is $k+n\Z$, where $k \in [0, n-1]$ is the number of blocks met by the subsequence. Since $1\le k\le s<n$, this sum is nonzero. The sequence is therefore zero-sum free, and its length, including the added basis terms, is
	\begin{equation}\label{girard:block-length}
		\sum_{i=1}^s(r_{j_i}+j_i)+\left(r+1-\sum_{i=1}^s r_{j_i}\right)
		=r+1+\sum_{i=1}^s j_i.
	\end{equation}
	
\smallskip
\noindent\emph{1. The density of exceptional ranks.}
We first show that $\mathsf d(G_r)\sim r+1$ as $r\to\infty$.

For each $i\in[1,\mathsf d(G_r)]$, let $b_i\in[0,n-1]$ be the unique integer satisfying $a_i=b_i+n\Z$. A vector $x=(x_1,\ldots,x_{\mathsf d(G_r)})\in\F_2^{\mathsf d(G_r)}$ selects the $i$th term of $S$ exactly when $x_i=1$. The linear map
$$
\varphi\colon\F_2^{\mathsf d(G_r)}\to\F_2^{r+1},\qquad
\varphi(x)=\sum_{i=1}^{\mathsf d(G_r)}x_i g_i
$$
gives the sum of the selected first coordinates. Set $V=\ker(\varphi)$. Since $\dim\varphi(\F_2^{\mathsf d(G_r)})\le r+1$, we obtain
$$
\dim V=\mathsf d(G_r)-\dim\varphi(\F_2^{\mathsf d(G_r)})
\ge\mathsf d(G_r)-(r+1).
$$

We set $M=\sum_{i=1}^{\mathsf d(G_r)}b_i$ and define
$$
\Phi\colon V\longrightarrow\F_2^M,\qquad
\Phi(x)=\bigl(\underbrace{x_1,\ldots,x_1}_{b_1},
\ldots,\underbrace{x_{\mathsf d(G_r)},\ldots,x_{\mathsf d(G_r)}}_{b_{\mathsf d(G_r)}}\bigr).
$$
The $i$th block has $b_i$ coordinates even when $x_i=0$, and is omitted if $b_i=0$. Index the coordinates of $\Phi(x)$ by $(i,h)$ with $1\le i\le\mathsf d(G_r)$ and $1\le h\le b_i$, so that $\Phi(x)_{i,h}=x_i$. Let $x,y\in V$ and $\alpha,\beta\in\F_2$. Since $V$ is a subspace, $\alpha x+\beta y$ belongs to $V$. For each coordinate,
$$
\Phi(\alpha x+\beta y)_{i,h}
=(\alpha x+\beta y)_i
=\alpha x_i+\beta y_i
=\bigl(\alpha\Phi(x)+\beta\Phi(y)\bigr)_{i,h}.
$$
Thus $\Phi$ is linear.

The Hamming weight $\mathsf{wt}(\Phi(x))$ equals $\sum_{i:\,x_i=1}b_i$. Indeed, if $x_i=1$ then the $i$th block contributes $b_i$ ones, and if $x_i=0$ it contributes none. For $x\in V\setminus\{0\}$, the selected subsequence is nontrivial and its first-coordinate sum is zero. Since $S$ is zero-sum free, its second-coordinate sum is nonzero, and in $\Z/n\Z$ we have
$$
0\ne\sum_{i:\,x_i=1}a_i
=\left(\sum_{i:\,x_i=1}b_i\right)+n\Z
=\mathsf{wt}(\Phi(x))+n\Z.
$$
Hence $\Phi(x)$ has weight not divisible by $n$ and cannot be zero. The linear map $\Phi$ is therefore injective. Thus $\dim\im(\Phi)=\dim V$, and zero is the only element of $\im(\Phi)$ with weight divisible by $n$.

Since $b_i\le n-1$, we have $M\le N=(n-1)\mathsf d(G_r)$. Appending $N-M$ zero coordinates to each element of $\im(\Phi)$ gives a subspace $C\subseteq\F_2^N$ of the same dimension and with the same weights. Recall that the PFR hypothesis of Lemma~\ref{girard:code-bound} follows from Lemmas~\ref{girard:torsion-pfr} and~\ref{girard:plunnecke}, as shown above. Since $N\ge\mathsf d(G_r)\ge r+1\ge2$, we can apply Lemma~\ref{girard:code-bound} to the affine subspace $C$ with ambient dimension $N$, modulus $n$ and residue zero to obtain
$$
\begin{aligned}
0\le\mathsf d(G_r)-(r+1)\le\dim V=\dim C
&\le\girardcodeD(N,n)\\
&\le c_n\frac{N}{\log N}
\le\frac{(n-1)c_n\mathsf d(G_r)}{\log\mathsf d(G_r)}.
\end{aligned}
$$
For the last inequality, we used $N\ge\mathsf d(G_r)$. Since $\mathsf d(G_r)\ge r+1$ tends to infinity and $n$ is fixed, division by $\mathsf d(G_r)$ gives
$$
\frac{\mathsf d(G_r)-(r+1)}{\mathsf d(G_r)}\longrightarrow0.
$$
Consequently $\mathsf d(G_r)\sim r+1$, and hence
$$
\frac{k_r}{r}=\frac{\mathsf d(G_r)}{r}-1\longrightarrow0.
$$

Suppose that  $r$ is exceptional and that the sequence  $S$, given in \eqref{given-sequence}, has an odd-order term $(0,a)$. Replace that term by $B_2(a)$ in the new direction. This block consists of three copies of $(e,b)$ with $2b=a$. If a nontrivial subsequence of the new sequence summed to zero, it would select an even number of these copies, since no other terms contribute in the new direction. If it selected none, it would be a nontrivial zero-sum subsequence of $S$. If it selected two, the two copies would sum to $(0,2b)=(0,a)$; replacing the pair by the original term would again give a nontrivial zero-sum subsequence of $S$. Thus the new sequence is zero-sum free of length  $\mathsf d (G_r)-1+3=\mathsf d (G_r)+2$, whence $k_{r+1}-k_r\ge1$.

For $x\ge 1$, set $m=\lfloor x\rfloor$. Each exceptional $r\in[1,m]$ contributes at least one to the following sum, and every other summand is nonnegative by \eqref{k_r-inequality}. Therefore
$$
  |\{r\in[1,m]: r\text{ is exceptional}\}|
  \le\sum_{r=1}^m(k_{r+1}-k_r)
  =k_{m+1}-k_1.
$$
Dividing by $m$, we obtain an upper bound of $(k_{m+1}-k_1)/m$ for the proportion of exceptional ranks in $[1,m]$. This bound tends to zero since $k_r/r\to0$. The complement of the exceptional ranks consequently has natural density one.

\smallskip
\noindent\emph{2. The case of a prime power.}
Let $n=p^b$ with $p$ being an odd prime. Suppose that  $S$ has an element of odd order, say
 $(g_{\mathsf d (G_r)} +a_{\mathsf d (G_r)})$ has odd order. We are going to prove that $\log_2 (r)$ is bounded above by a constant depending on $n$. Thus, for every fixed $n$, $S$ has no elements of odd order for all sufficiently large $r$.

Since $(g_{\mathsf d (G_r)} +a_{\mathsf d (G_r)})$ has odd order, it follows that   $g_{\mathsf d (G_r)}=0$, whence $a_{\mathsf d (G_r)}\ne0$ because $S$ is zero-sum free. We first derive an inequality that does not require $r$ to be large. We define a map
$$
  A'\colon\F_2^{\mathsf d (G_r)-1}\longrightarrow\F_2^{r+1},\qquad
  A'(x)=\sum_{i=1}^{\mathsf d (G_r)-1}x_i g_i \quad \text{with} \ x=(x_1,\ldots,x_{\mathsf d(G_r)-1})\in\F_2^{\mathsf d(G_r)-1} \,.
$$
We set $W=\ker (A')$ and obtain that $\dim W\ge \mathsf d (G_r)-1-(r+1)$. For $i \in [1,\mathsf d(G_r)-1]$, let $\lambda_i\colon\F_2^{\mathsf d (G_r)-1}\to \{0,1\}\subset\Z$ send each vector to the integer representative of its $i$th coordinate. For $x\in W$, let
$$
  F(x)=\sum_{i=1}^{\mathsf d (G_r)-1}\lambda_i(x)a_i\in  \Z / n\Z.
$$
For $x\ne0$ in $W$, the selected terms form a nontrivial subsequence whose first-coordinate sum is zero. Omitting the last term shows $F(x)\ne0$; including it also gives first-coordinate sum zero, so $F(x)+a_{\mathsf d (G_r)}\ne0$. For $x=0$, we have $F(0)=0$ and $F(0)+a_{\mathsf d (G_r)}=a_{\mathsf d (G_r)}\ne0$.

 We use the group algebra $\F_p[\Z / n\Z]$ and all notation, introduced in \eqref{eq:augmentation-ideal}, \eqref{vanishing-of-the-augmentation-ideal}, and \eqref{basic-notation}. The preceding observations show that
$$
  [X^0]\bigl((1-X^{a_{\mathsf d (G_r)}})X^{F(x)}\bigr)
  =[X^0]\bigl(X^{F(x)}-X^{F(x)+a_{\mathsf d (G_r)}}\bigr)
  =\mathbf1_{\{0\}}(x)\qquad(x\in W).
$$
Both terms have exponents different from zero when $x\ne0$, while for $x=0$ only the first term has exponent zero.

We use expansion \eqref{growth:subset-expansion} with $G=\Z / n \Z$, $h_i=a_i$, $q=p$, and $\ell=\mathsf d (G_r)-1$.  The $i$th factor below is $1$ if $x_i=0$ and $X^{a_i}$ if $x_i=1$.
$$
\begin{aligned}
  X^{F(x)}
    &=\prod_{i=1}^{\mathsf d (G_r)-1}\bigl(1+(X^{a_i}-1)\bigl(\lambda_i(x)+p\Z\bigr)\bigr)\\
    &=\sum_{J\subseteq[1,\mathsf d (G_r)-1]}
       \left(\prod_{i\in J}(X^{a_i}-1)\right)
       \left(\prod_{i\in J}\bigl(\lambda_i(x)+p\Z\bigr)\right).
\end{aligned}
$$
In the expansion, $J$ records the positions at which the second summand of the corresponding factor was taken. The products for $J=\emptyset$ equal $1$.

Let $I_{\Z / n \Z} \subseteq \F_p[\Z / n \Z]$ be the augmentation ideal. Then
$$
  (1-X^{a_{\mathsf d (G_r)}})\prod_{i\in J}(X^{a_i}-1)\in I_{\Z/ n \Z}^{|J|+1}.
$$
Since $\Z / n \Z$ is a $p$-group, \eqref{vanishing-of-the-augmentation-ideal} and $\mathsf D(\Z/n \Z)=n$ give $I_{\Z / n \Z}^n=\{0\}$. Hence all the displayed coefficients with $|J|\ge n-1$ vanish. For the remaining $J$, set
$$
  c_J=[X^0]\left((1-X^{a_{\mathsf d (G_r)}})\prod_{i\in J}(X^{a_i}-1)\right)
  \in\F_p .
$$
Multiplying the expansion by $1-X^{a_{\mathsf d (G_r)}}$ and taking the coefficient of $X^0$ gives
$$
  \mathbf1_{\{0\}}(x)
  =\sum_{\substack{J\subseteq[1,\mathsf d (G_r)-1]\\ |J|\le n-2}}c_J\prod_{i\in J}\bigl(\lambda_i(x)+p\Z\bigr).
$$
Thus the indicator on $W$ is obtained by evaluating a polynomial of degree at most $n-2$ in the $\mathsf d (G_r)-1$ coordinate variables.

For $w\in W$, addition by $w$ preserves the $i$th bit if $w_i=0$ and interchanges $0$ and $1$ if $w_i=1$. In terms of the integer representatives, $\lambda_i(x+w)$ is therefore $\lambda_i(x)$ or $1-\lambda_i(x)$, respectively, for every $x\in W$. This verifies the translation hypothesis of Lemma~\ref{prelim:cosets}. Apply that lemma with ambient space $\F_2^{\mathsf d (G_r)-1}$, the subspaces $\{0\}\le W$, coefficient field $\F_p$, $\ell=\mathsf d (G_r)-1$, and degree bound $t=n-2$. Using the lower bound for $\dim W$, we obtain
\begin{equation}\label{girard:marked}
  2^{\mathsf d (G_r)-1-(r+1)}\le 2^{\dim W}\le\sum_{j=0}^{n-2}\binom{\mathsf d (G_r)-1}{j}.
\end{equation}

For $j \in [0, n-2]$, we have $\binom{\mathsf d (G_r)-1}{j}\le {\mathsf d (G_r)}^j\le {\mathsf d (G_r)}^{n-2}$. Taking base-two logarithms of \eqref{girard:marked} now gives
$$
  \mathsf d (G_r)-(r+1) \le 1+\log_2(n-1)+(n-2)\log_2 \mathsf d (G_r).
$$
By the asymptotic estimate established in Part~1, $\mathsf d (G_r)\sim r+1$, so $\mathsf d (G_r)\le 2(r+1)$ for all sufficiently large $r$. Thus, under the odd-order assumption, we obtain the upper bound
$$
  \mathsf d (G_r)-(r+1) \le (n-2)\log_2 (r+1)+1+\log_2(n-1)+(n-2).
$$

Independently of that assumption, we obtain a lower bound from the block construction. Take $n-1$ copies of $B_j(1)$, with
$$
  j=\left\lfloor\log_2\frac{r+1}{n-1}\right\rfloor .
$$
For sufficiently large $r$ this is at least $1$. It is the largest integer satisfying the sufficient dimension condition $(n-1)2^j\le r+1$: the actual dimension required by the blocks is $(n-1)r_j\le(n-1)2^j$, whereas their resulting length $r+1+(n-1)j$ increases with $j$. By \eqref{girard:block-length} and the maximality of $\mathsf d (G_r)$, we obtain
$$
  \mathsf d (G_r)- (r+1) \ge(n-1)j
  \ge (n-1)\log_2 (r+1)-(n-1)\bigl(\log_2(n-1)+1\bigr).
$$
The last inequality uses $\lfloor y\rfloor\ge y-1$. If both bounds hold, subtracting $(n-2)\log_2 (r+1)$ shows that $\log_2 (r+1)$ is at most a constant depending only on $n$.

\smallskip
\noindent\emph{3. The case $n=3$.} We proceed in two steps.

3.1 Suppose  that the sequence $S$, given in \eqref{given-sequence}, has an element of   odd order.  By \eqref{basic-1}, we have $\mathsf d (G_r)\ge r+5$. Put $j=\mathsf d (G_r)-(r+1)-1$, so $j\ge3$ and $\mathsf d (G_r)=r+1+j+1$. Since $n=3$ is a prime power, \eqref{girard:marked} applies and gives
$$
  2^j=2^{\mathsf d (G_r)-1-(r+1)}\le \binom{\mathsf d (G_r)-1}{0}+\binom{\mathsf d (G_r)-1}{1}=\mathsf d (G_r).
$$
Consequently $r_j=2^j-1-j\le \mathsf d (G_r)-1-j=r+1$.

Consider the two blocks $B_{j-1}(1)$ and $B_3(1)$. Their first-coordinate spaces have dimensions $r_{j-1}$ and $r_3=2^3-1-3=4$. Assume to the contrary that  $r+1\ge r_{j-1}+4$. Then the two blocks fit in $\F_2^{r+1}$. Since $2<n=3$, the block construction gives a zero-sum free sequence whose length, by \eqref{girard:block-length}, is
$$
  r+1+(j-1)+3=r+1+j+2=\mathsf d (G_r)+1,
$$
a contradiction. Thus, it follows that $r_j\le r+1<r_{j-1}+4$. Subtracting $r_{j-1}$ gives
$$
  2^{j-1}-1
  =(2^j-1-j)-(2^{j-1}-1-(j-1))
  =r_j-r_{j-1}<4.
$$
Since $j\ge3$, this implies that $j=3$. Therefore $r_2=2^2-1-2=1$, whence $4=r_3\le r+1<r_2+4=5$, and so we get  $r=3$.

\smallskip
3.2 Let $r=3$, whence $G_3 = C_2^3 \oplus C_6$. Take $B_3(1)$ and $B_1(1)$, whose first-coordinate spaces have dimensions $4$ and $0$. They use all $4$ dimensions, so there are no added basis terms. The construction gives a zero-sum free sequence of length $4+3+1=8$. Its block $B_1(1)$ consists of the single term $0+(1+n\Z)$ of order $3$.

By \cite[Theorem 3.9.6]{Ge-Gr-Zh26a}, we have  $\mathsf D(G_3)=\mathsf D^*(G_3) = 9$.
Hence the sequence of length $8$ is extremal. Adjoining the negative of its sum gives a minimal zero-sum sequence of length $9$, which contains $(0,1)$, as required.
\end{proof}
\endgroup

\providecommand{\bysame}{\leavevmode\hbox to3em{\hrulefill}\thinspace}
\providecommand{\MR}{\relax\ifhmode\unskip\space\fi MR }
\providecommand{\MRhref}[2]{%
  \href{http://www.ams.org/mathscinet-getitem?mr=#1}{#2}
}
\providecommand{\href}[2]{#2}

\end{document}